\documentclass[a4paper,english,10pt,oneside]{article} %zibreporttrue

\usepackage{ifthen}
\newif\ifzibreport
\newif\ifarxiv

\arxivtrue

\zibreporttrue
\ifzibreport
\usepackage{amsfonts,amssymb,mathrsfs}
\else
\usepackage{IJOC_style/eqndefns-left} % For checking the display equation width and equation environment definitions %
\OneAndAHalfSpacedXII % Current default line spacing
\RequirePackage{tgtermes}
\RequirePackage{newtxtext}
\RequirePackage{newtxmath}
\RequirePackage{bm}
\RequirePackage{endnotes}
\fi

\usepackage{algorithm}
\usepackage{algpseudocode}
\usepackage{tikz}
\usepackage{longtable}
\usepackage{booktabs}   % optional, but keeps rules consistent

\usepackage{ifthen}
\usepackage{adjustbox}
\usepackage{xspace}
\usepackage{booktabs}
\usepackage{tabularx}
\usepackage{amsmath}
\usepackage{subcaption}
\usepackage[textsize=scriptsize, textwidth=3.3cm]{todonotes}
\usepackage{hyperref}
\usepackage{tikz}
\usetikzlibrary{matrix,calc,positioning}
\usetikzlibrary{petri}	% for tokens (dots)
\usetikzlibrary{shapes.misc}
\tikzset{cross/.style={cross out, draw=black, minimum size=2*(#1-\pgflinewidth), inner sep=0pt, outer sep=0pt},
	cross/.default={1pt}}
\usepackage{pgfplots}
\pgfplotsset{compat=1.16}

\ifzibreport
\NewDocumentEnvironment{bMatrix}{m}{\begin{bmatrix}}{\end{bmatrix}}
\let\ABSTRACT\abstract
\let\FIGURE\relax
\let\TABLE\relax
\newcommand{\ACKNOWLEDGMENT}[1]{%
	\begin{acknowledgement}
		#1
	\end{acknowledgement}%
}
\else
\fi

\newcommand{\myorcidlink}[1]{\,\href{https://orcid.org/#1}{\raisebox{-0.45ex}{\includegraphics[width=1.8ex]{orcid}}}}
\newcommand{\myurl}[1]{\textsf{\footnotesize \url{#1}}\xspace}

\newcommand{\allcaps}[1]{\protect\scalebox{0.93}{#1}}

\newcommand{\name}[1]{\mbox{#1}\xspace}
\newcommand{\nameCaps}[1]{\name{\allcaps{#1}}}

\newcommand{\MPI}{\nameCaps{MPI}}

\newcommand{\T}{^{T}}

\newcommand{\subjectto}{\mbox{subject to}\xspace}

\newcommand{\N}{\mathbb{N}}
\newcommand{\R}{\mathbb{R}}

\definecolor{seagreen}{rgb}{0.18,0.74,0.56}
\definecolor{darkgreen}{rgb}{0.0,0.45,0.00}
\definecolor{navyblue}{rgb}{0.0,0.0,0.5}
\definecolor{steelblue}{rgb}{0.27,0.51,0.71}
\definecolor{siennabrown}{rgb}{0.63,0.32,0.18}
\definecolor{firebrickred}{rgb}{0.69,0.13,0.13}
\definecolor{gray75}{rgb}{0.75,0.75,0.75}
\definecolor{orange}{rgb}{.843,0.671,0.078}
\definecolor{gold}{rgb}{1.0,0.84,0.0}
\definecolor{scipyellow}{HTML}{FFFFD6}
\definecolor{soplexred}{HTML}{FFD8D8}
\definecolor{zimplgreen}{HTML}{D8FFD8}
\definecolor{ugblue}{HTML}{CFEFFF}
\definecolor{gcgorange}{HTML}{FFDAB9}

\definecolor{mDarkBrown}{HTML}{604c38}
\definecolor{mDarkTeal}{HTML}{23373b}
\definecolor{mLightBrown}{HTML}{EB811B}
\definecolor{mLightblue}{HTML}{14B03D}

\definecolor{c0}{HTML}{000060}
\definecolor{c1}{HTML}{0000FF}
\definecolor{c2}{HTML}{36648B}
\definecolor{c3}{HTML}{4682B4}
\definecolor{c4}{HTML}{5CACEE}
\definecolor{c5}{HTML}{FF0000}
\definecolor{c6}{HTML}{008888}
\definecolor{c7}{HTML}{00DD99}
\definecolor{c8}{HTML}{527B10}
\definecolor{c9}{HTML}{7BC618}
\definecolor{c10}{HTML}{8DD8F8}
\definecolor{background}{HTML}{FFFFFF}

\definecolor{cp1}{HTML}{F2AF29}
\definecolor{cp2}{HTML}{05668D}
\definecolor{cp3}{HTML}{02C39A}
\definecolor{cp4}{HTML}{4F345A}
\definecolor{cp5}{HTML}{F2B5D4}
\definecolor{cp6}{HTML}{DA2C38}

\usepackage{natbib}
\bibpunct[, ]{(}{)}{,}{a}{}{,}%

\ifzibreport\else
\EquationsNumberedThrough    % Default: (1), (2), ...
\TheoremsNumberedThrough     % Preferred (Theorem 1, Lemma 1, Theorem 2)
\ECRepeatTheorems  %  
\MANUSCRIPTNO{IJOC-0001-2024.00}
\fi

\newcommand{\TheTitle}{Distributed Parallel Structure-Aware Presolving for Arrowhead Linear Programs}
\newcommand{\TheAuthors}{N.--C. Kempke, S. Maher, D. Rehfeldt, A. Gleixner, T. Koch, S. Uslu}

\hypersetup{%
	pdftitle={\TheTitle},
	pdfauthor={\TheAuthors}
}

\ifzibreport
	\usepackage{./zibtitlepage}
	
    \providecommand{\Description}[1]{}
	\newenvironment{acknowledgement}{\section*{Acknowledgements}}
\fi

\begin{document}
%%%%%%%%%%%%%%%%

\ifzibreport
	\ZTPTitle{\TheTitle}
	\title{\TheTitle}
	
	\ZTPAuthor{
		\ZTPHasOrcid{Nils--Christian Kempke}{0000-0003-4492-9818},
		\ZTPHasOrcid{Stephen J. Maher}{0000-0003-3773-6882},
		\ZTPHasOrcid{Daniel Rehfeldt}{0000-0002-2877-074X},
		\ZTPHasOrcid{Ambros Gleixner}{0000-0003-0391-5903},
		\ZTPHasOrcid{Thorsten Koch}{0000-0002-1967-0077},
		Svenja Uslu
	}
	\author{
		\ZTPHasOrcid{Nils--Christian Kempke}{0000-0003-4492-9818},\and\
		\ZTPHasOrcid{Stephen J. Maher}{0000-0003-3773-6882},\and\
		\ZTPHasOrcid{Daniel Rehfeldt}{0000-0002-2877-074X},\and\
		\ZTPHasOrcid{Ambros Gleixner}{0000-0003-0391-5903},\and\
		\ZTPHasOrcid{Thorsten Koch}{0000-0002-1967-0077}, \and\
		Svenja Uslu
    }
	
	\ZTPInfo{Preprint}
	\ZTPNumber{26-01}
	\ZTPMonth{March}
	\ZTPYear{2026}
	
	\date{\normalsize Match 03, 2026}
	\ifarxiv\else
		\zibtitlepage
	\fi
	\maketitle
\fi

% TODO: Abstract <= 300 words!
\ABSTRACT{%
We present a structure-aware parallel presolve framework for arrowhead linear programs (AHLPs) in high-performance computing (HPC) environments, integrated into the parallel interior point solver PIPS-IPM++. Large-scale LPs arising from automated model generation frequently contain redundancies and numerical pathologies that necessitate effective presolve, yet existing presolve techniques are primarily serial or structure-agnostic. Our presolve works on the distributed data of PIPS-IPM++, communicates only information on the linking parts of the problem, and preserves the arrowhead structure. Already on a single compute node with 64 processes, it outperforms PaPILO by a factor of 19 and Gurobi's presolve by a factor of 6 (shifted geometric mean), while removing a few percentage points fewer nonzeros. On large instances with small linking parts, multiple compute nodes make it 20 to 115 times faster than Gurobi's presolve. Presolve raises the number of instances solved by PIPS-IPM++ from 37 to 48. External presolve works with the modeler's annotation carried over, but it is slower, especially on large models, and breaks the distributed workflow. With automatic re-annotation, PIPS-IPM++ solves at most 21 instances.
}%

\ifzibreport\else
%\FUNDING{This research was supported by [grant number, funding agency].}
%Supplemental Material:
%Data Ethics & Reproducibility Note:
%\KEYWORDS{Stochastic programming, Decision support,Uncertainty, Disaster response, Optimization}
%\HISTORY{Received: Month DD, YYYY; Accepted: Month DD, YYYY; Published Online: Month DD, YYYY}
\fi

% TODO: properly follow instructions on extended proceedings paper as stated here: https://pubsonline.informs.org/page/ijoc/submission-guidelines
% 25 pages limit + 10 pages appendices
% TODO: check keyword requirements
\maketitle
%%%%%%%%%%%%%%%%%%%%%%%%%%%%%%%%%%%%%%%%%%%%%%%%%%%%%%%%%%%%%%%%%%%%%%

% Text of your paper here
\section{Introduction}\label{sec:intro}

Presolve is an indispensable component of modern linear programming (LP) solvers. Its primary purpose is not only to reduce solution time, but also to improve numerical robustness and algorithmic stability by tightening formulations, eliminating redundancies, and correcting numerical pathologies before the actual solution process begins. In large-scale LPs, these effects are often decisive: poor conditioning, redundant constraints, or rank deficiencies can severely degrade or even prevent the convergence of simplex, interior-point methods (IPMs), or first-order methods (FOMs). Importantly, such deficiencies are rarely the result of modeling errors. In practice, large-scale LPs are typically generated automatically using high-level modeling tools and scenario-based expansion, where redundancies and linear dependencies arise naturally as a by-product of modular and expressive model construction. Consequently, presolve has become a standard preprocessing step in all state-of-the-art LP solvers.

At the same time, the size and structural complexity of LPs arising in practice have grown dramatically. Many contemporary applications, particularly in energy systems, stochastic programming, and multi-stage planning, give rise to problems with millions of variables and constraints and highly structured constraint matrices. In such settings, solver scalability increasingly depends on the ability to exploit problem structure and parallel hardware effectively. While considerable progress has been made in developing structure-exploiting and distributed solution algorithms, presolve has largely remained a serial, structure-agnostic preprocessing step. As a result, presolve itself can become a computational bottleneck and, more critically, may interfere with or even destroy the very structure that specialized solvers rely on for scalability.

A recurring structural pattern in LPs is the \emph{arrowhead} or \emph{doubly bordered block-diagonal form}, as illustrated in Figure~\ref{fig:remix_blockstructure}. Arrowhead LPs (AHLPs) contain \emph{linking variables} (connecting diagonal blocks vertically) and \emph{linking constraints} (connecting blocks horizontally). This structure generalizes both primal- and dual block-angular problems.
\ifzibreport
\begin{figure}
	\centering
	\includegraphics[width=0.4\linewidth]{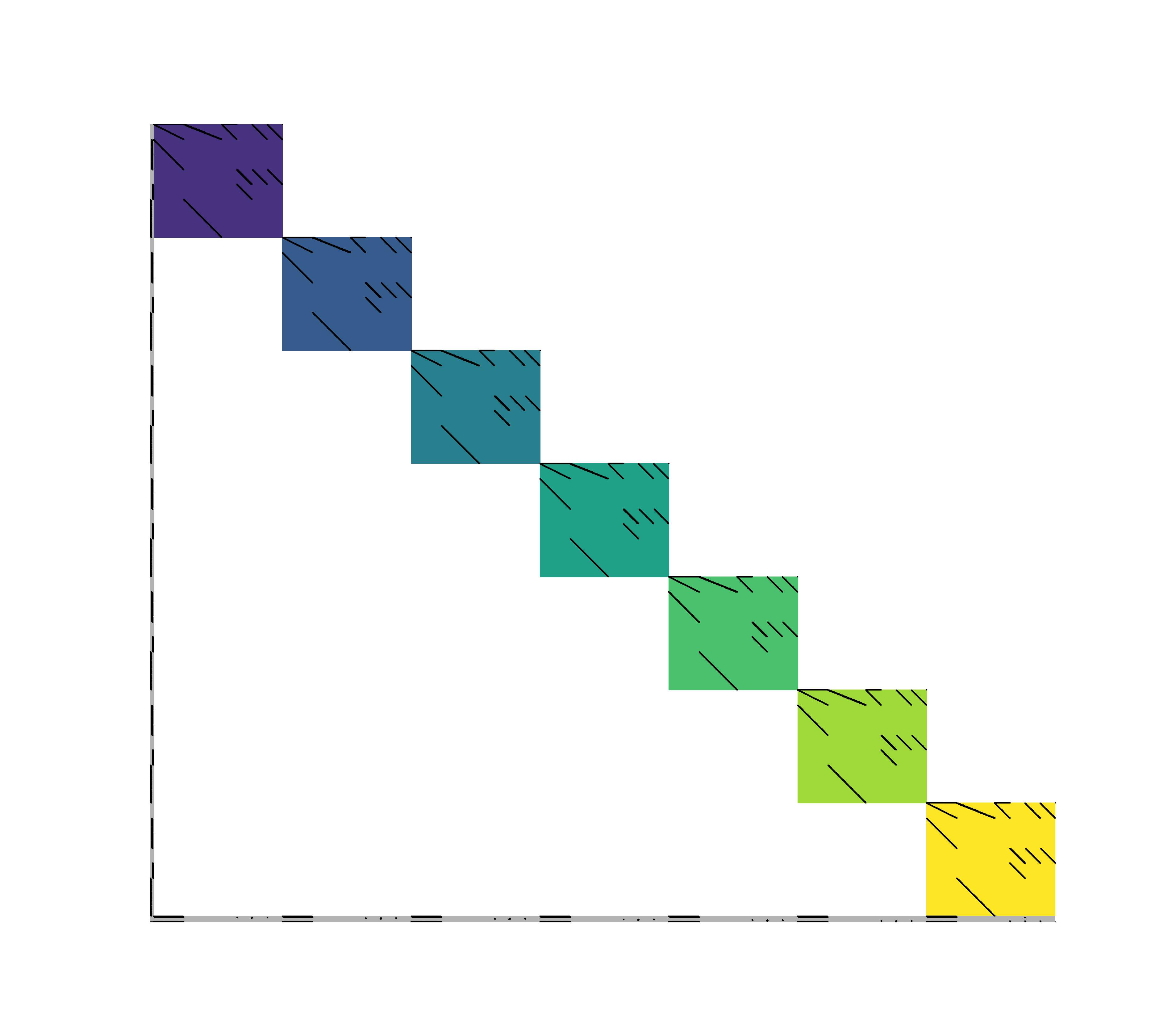}
	\caption{Constraint matrix with the arrowhead structure of a real-world energy system model.}
	\label{fig:remix_blockstructure}
\end{figure}
\else
\begin{figure}
	\FIGURE
	{\includegraphics[width=0.4\linewidth]{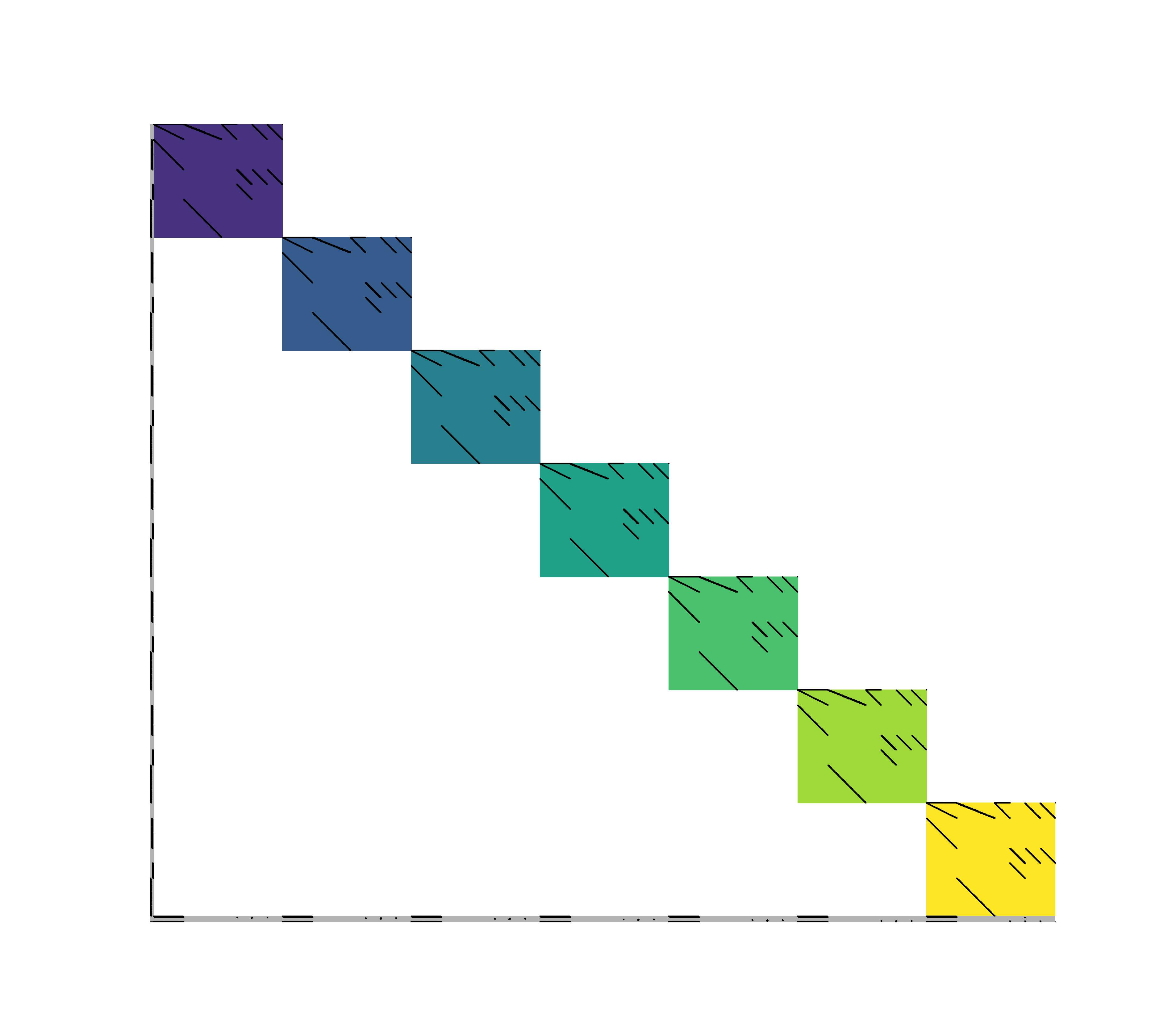}}
	{Constraint matrix with the arrowhead structure of a real-world energy system model.\label{fig:remix_blockstructure}}
	{}
\end{figure}
\fi
AHLPs arise in a wide range of practical applications. Energy system models (ESMs), such as electricity market models with dispatch decisions~\cite{Rehfeldt2019_PIPSIPMpp}, renewable expansion planning~\cite{REMIX_2017,PyPSAEur}, and large-scale (stochastic) economic (re-)dispatch models, can be structured by grouping variables and constraints along spatial or temporal dimensions, resulting in an arrowhead LP. Multi-stage LPs, where decisions at each stage depend on preceding stages, and multi-stage stochastic LPs, covering applications such as asset–liability management, supply network design, revenue management, and portfolio optimization, exhibit arrowhead structure~\cite{Castro_MultiStage,Colombo2011,Steinbach2001}. Staircase LPs used in production scheduling, inventory management, transportation, and multistage system design show local linking structure connecting consecutive diagonal blocks~\cite{Fourer1982_StaircaseLPs,Wittrock1985_DualNestedOfStaircaseLP}. Band-diagonal matrices arising in distribution planning can also be reformulated to expose arrowhead forms~\cite{OOPS_DistributionPlanning}.

The arrowhead structure also arises in mixed-integer ESMs~\cite{WetzelGreenEnergyCarriers2023} and in nonlinear optimization, e.g., in portfolio optimization~\cite{GondzioOOPS_NonlinPortfolioOptimization}, dynamic optimization~\cite{Word2014_EfficientParallelIPMDAT}, and multi-stage nonlinear programs~\cite{MadNLPGPU}.

Arrowhead structures can often be exploited algorithmically with specialized solution methods. Structure-exploiting algorithms promise to overcome scaling issues and push current computational limits by leveraging high-performance computing (HPC), where general LP solution techniques often scale poorly. Solution techniques for AHLPs include specialized simplex methods~\cite{Fourer1982_StaircaseLPs,Friedlander1990_StaircaseLPSimplex}, Dantzig-Wolfe decomposition~\cite{Wittrock1985_DualNestedOfStaircaseLP}, Benders decomposition~\cite{Meersman2023_NestedBendersStochastic,Zhang2024_BendersForLargeStochasticLP}, Lagrangian decomposition~\cite{kim19}, and interior-point methods (IPMs)~\cite{BorderedBlockDiagIPM1996,LubinDualDecomp2013}. IPM solvers that explicitly exploit this structure include BlockIp~\cite{Castro2007_BlockIP}, OOPS~\cite{Gondzio2008_OOPS2}, PIPS-IPM~\cite{Petra2014_PIPSaugmented}, and MadNLP~\cite{MadNLPGPU}. While these solvers routinely operate in distributed parallel environments, presolve is typically applied as a monolithic preprocessing step prior to the parallel solution phase, often in a simplified form to avoid destroying the arrowhead structure.

All solution techniques, whether structure-exploiting or not, run more robustly and efficiently on tight formulations free of redundancies and numerical pathologies. Presolve reformulates a given model by applying a sequence of reduction techniques, generating an equivalent but simpler and smaller \emph{presolved} model. This presolved model is then passed to the solution algorithm, solved, and in the subsequent \emph{postsolve} step, a solution to the original problem is obtained by reverting the presolve reductions. Presolve for LP has a long history~\cite{AndersenAndersenPresolve1995,GondzioPresolveLP1997,AchterbergBixbyGuetal2016,TwoRowPresolve_Gemander2020}, and today, all commercial LP solvers and most academic LP solvers apply presolve before starting the solution algorithm.

Combining presolve with structured LPs and specialized solvers is challenging: presolve must preserve enough of the underlying problem structure for the solver to work efficiently. While general-purpose presolve libraries exist, most notably the academic parallel presolve library PaPILO~\cite{Gleixner2023_Papilo}, they do not keep the assignment of constraints and variables to blocks, so the structure has to be recovered after presolve. Additionally, these libraries fail to efficiently leverage the arrowhead structure. This creates a clear need for presolve techniques that are both parallel and structure-aware. In this paper, we present our implementation of specialized presolve routines within the distributed parallel IPM PIPS-IPM++. Our presolve efficiently distributes reductions across compute nodes, applies presolve routines in a distributed, parallel manner, communicates only information on the linking parts of the problem, and preserves the structure of AHLPs.

\subsection{Contribution}

We make three contributions: a structure-aware distributed parallel presolve and postsolve framework for AHLPs, an experimental evaluation demonstrating its efficiency and scalability on a set of large-scale AHLPs, and a comparison with the parallel presolve library PaPILO and Gurobi's presolve implementation, both as standalone presolve and as external presolve for PIPS-IPM++. Already on a single compute node with 64 processes, our presolve outperforms Gurobi and PaPILO by factors of 6 and 19, respectively, and on multiple compute nodes it is up to 115 times faster than Gurobi's presolve on large instances. An external presolve requires recovering the block structure, which an automatic annotation cannot do well. Even with the modeler's annotation carried over, it is not faster overall, costs much more time on large instances, and breaks the distributed pipeline of model generation and solution.

\paragraph{Previous proceedings publication} Early experiments on a first version of our presolve were published in the conference proceedings~\cite{Kempke19}. Since then, most of our implementation has been rewritten and redesigned. This paper extends \cite{Kempke19} by a fully implemented linear presolve framework, in-depth implementation and design details and a substantially extended computational study on a large set of AHLPs.

\subsection{Message Passing Interface}\label{subsec:mpi}

Our parallel implementation uses the Message Passing Interface (MPI)~\cite{LydonEtAl1994_MPIStandard}, the de facto standard for communication between processes on distributed-memory systems. We refer to the processes participating in MPI communication as \emph{MPI processes} (or simply \emph{processes}).

MPI provides collective operations such as \texttt{Reduce} and \texttt{Allreduce}, which apply a reduction operator (e.g., summation) to data contributed by all processes and return the result to a designated root process or to all processes, respectively. We use the blocking variants, which return once the calling process’s participation is complete. Note that \texttt{Allreduce} is effectively synchronizing, since its result depends on contributions from every process. MPI supports custom reduction operators, which we exploit in the parallel constraint detection (Section~\ref{subsec:parallel_rows}).

\subsection{Outline and Notation} 

\paragraph{Outline} Section~\ref{sec:presolve} introduces presolving and our structure-aware presolving framework, Section~\ref{sec:algos} describes the distributed implementation, highlighting synchronization and communication, Section~\ref{sec:results} evaluates it, and Section~\ref{sec:conclusion} concludes.

\paragraph{Notation}
Zero entries in large block matrices are sometimes omitted. When shown explicitly, they are denoted by $0$ with the size inferred from the context. For a matrix $A\in \R^{m\times n}$, we denote by $A_{i.}$ its $i$-th row and by $A_{.j}$ its $j$-th column.

\section{Structure-aware presolving techniques}\label{sec:presolve}

We consider linear optimization problems of the form
\begin{equation}\label{eq:pips_lp}
\min_x \; c^T x \quad \subjectto \quad A x = b, \quad d \le C x \le f, \quad l \le x \le u,
\end{equation}
with $x, c, l, u\in\R^n$, $A\in\R^{m_A \times n}$, $b\in\R^{m_A}$, $C\in\R^{m_C \times n}$, and $d, f\in\R^{m_C}$, where the entries of $d$, $f$, $l$, and $u$ may be infinite. The dual problem of Equation~(\ref{eq:pips_lp}) is
\begin{equation} \label{eq:pips_lp_dual}
\begin{aligned}
\max_{y,z^+,z^-,\gamma, \phi}\quad & b^T y + d^T z^+ - f^T z^- + l^T \gamma - u^T \phi \\
\subjectto \quad & A^T y + C^T (z^+ - z^-) + \gamma - \phi = c, \quad z^+, z^-, \gamma, \phi \ge 0,
\end{aligned}
\end{equation}
where $y \in \R^{m_A}$ (free), $z^+$, $z^- \in \R^{m_C}$, and $\gamma$, $\phi \in \R^n$ are the duals of the equality constraints, the inequality constraints, and the variable bounds, and $z^+_i$, $z^-_i$, $\gamma_i$, and $\phi_i$ are fixed at zero for $d_i=-\infty$, $f_i=\infty$, $l_i=-\infty$, and $u_i=\infty$, respectively.

A variety of presolving techniques exist \cite{AndersenAndersenPresolve1995,GondzioPresolveLP1997,AchterbergBixbyGuetal2016}. In our framework, these techniques are implemented as \emph{presolvers}. Each presolver applies one reduction or a combination of similar reductions. To do so, a presolver usually iterates over the problem's constraints, variables, non-zeros, or a subset of these. Presolvers are applied iteratively in \emph{rounds} until no further reductions are possible or predefined working limits are reached. Working limits are imposed to restrict the time spent in presolve. Within our solver, we impose a maximum number of presolving rounds, and a new round is only started if a sufficient number of reductions were applied in the previous round relative to the problem size. This prevents costly scans of the entire problem when only minimal additional reductions are expected.

\subsection{Presolvers in PIPS-IPM++}

% Define the common table rows
\newcommand{\presolverrows}{%
	\toprule
Presolver & Description \\
\midrule
Aggregation     & Aggregates two variables using doubleton constraints (\cite{GondzioPresolveLP1997} 2.2) \\
BoundTightening & One-constraint activity-based bound tightening (\cite{AchterbergBixbyGuetal2016} 3.2) \\
VarsFixation    & Fixes variables with (nearly) equal bounds (\cite{AchterbergBixbyGuetal2016} 4.1) \\ 
DualTightening  & Dual variable fixing and bound tightening (\cite{AchterbergBixbyGuetal2016} 4.4) \\
EmptyVar        & Removes empty variables (\cite{AndersenAndersenPresolve1995} 3.1 (ii)) \\
ForcingConstr   & Fixes variables in forcing constraints (\cite{AndersenAndersenPresolve1995} 3.3 (x)) \\
LinDependencies & Detects linear dependencies in the equality matrix (Section~\ref{sec:lindep}, \cite{GondzioPresolveLP1997} 2.1 and 3.1) \\ 
ParallelVars    & Detects and merges parallel variables (\cite{AndersenAndersenPresolve1995} 3.4) \\
ParallelConstrs & Detects parallel or nearly parallel constraints (\cite{AchterbergBixbyGuetal2016} 5.2) \\
Permutation     & Improves arrowhead structure by permuting variables and constraints (Section~\ref{sec:permutation}) \\ 
RedundantExpr   & Removes redundant expressions (\cite{Achterberg2013}, 4.4 (i)) \\
RedundantConstr & Removes redundant (w.r.t. activity) constraints (\cite{AchterbergBixbyGuetal2016} 3.1) \\
SingletonVar    & Substitutes singleton variables (\cite{AndersenAndersenPresolve1995} 3.2 (vi--viii)) \\ 
SingletonConstr & Transforms singleton constraints into variable bounds (\cite{AndersenAndersenPresolve1995} 3.1 (v))\\
TinyEntries     & Removes small entries in the constraint matrix (\cite{AchterbergBixbyGuetal2016} 3.1) \\
\bottomrule
}

\ifzibreport
\begin{table}
	\begin{tabularx}{\textwidth}{lX}
		\presolverrows
	\end{tabularx}
	\caption{Presolving techniques in PIPS-IPM++.}
	\label{tab:presolvers}
\end{table}
\else
\begin{table}
	\TABLE
	{Presolving techniques in PIPS-IPM++.\label{tab:presolvers}}
	{\begin{tabular}{ll}
			\presolverrows
	\end{tabular}}
	{}
\end{table}
\fi

The presolvers implemented in PIPS-IPM++ are shown in Table~\ref{tab:presolvers}. The presolvers EmptyVar, ForcingConstr, RedundantConstr and TinyEntries are actually part of PIPS-IPM++'s ModelCleanup presolver, but for ease of presentation of the individual techniques, we disaggregated them in Table~\ref{tab:presolvers}. The presolver RedundantExpr runs as part of the DualTightening presolver as it uses the same detection mechanism. The presolvers we implemented are not an exhaustive list of all available presolving techniques. Rather, we implemented all presolvers that seemed general and widely adopted. Additionally, we selected presolvers that were expected to perform reductions on the models of interest, which primarily come from energy modeling contexts. We omit a detailed explanation of each presolver's reduction technique, except for the LinDependencies and Permutation presolvers. We instead provide references to the respective papers in the ``Description'' column of Table~\ref{tab:presolvers}.

All of the presolvers are called in each round of the PIPS-IPM++ presolve routine. The exceptions are the Tiny Entries presolver, which is called once at the beginning of presolve, and the LinDependencies presolver, which is called once at the end of presolve. The presolvers Permutation and LinDependencies are specifically tailored to AHLPs and will be described in more detail in Section~\ref{sec:permutation} and Section~\ref{sec:lindep}.

\subsection{Distributed Arrowhead presolving}

Within PIPS-IPM++, the LP in Equation~(\ref{eq:pips_lp}) is of \emph{arrowhead} (primal-dual block-angular) form:
\begin{equation}\label{eq:pips_blocklp}
\begin{alignedat}{4}
    \min \quad {{c_0^T x_0} }~~+~~			& { c_1^T x_1 }~~+~~\cdots		&~~+~~c_N^T x_N		& \\[-6pt]
    \subjectto \quad {A_0 x_0} ~~~~~~~ 			&   								 	 		&										& = {b_0} \\[-6pt]
    {d_0} \leq{C_0  {x_0}} ~~~~~~~ 	&  									  			&										&\leq{f_0}  \\[-6pt]
    {A_1 {x_0} }~~+~~     						&  B_1 x_1 		 		        &        								&= b_1  \\[-6pt]
    { d_1}\leq{  C_1  {x_0}}~~+~~& D_1 x_1 			            &        								&\leq f_1  \\[-6pt]
    { \vdots}~~~~~~~~~~											&~~~~~~~~~~~~~~\ddots     &   	  														&~\vdots  \\[-6pt]
    { 	A_N{x_0}} ~~+~~      					&        							            &~~~~~~ B_N x_N	  		&= b_N  \\[-6pt]
    { 	d_N	}\leq{C_N {x_0}}~~+~~&				                    	       	&~~~~~~ D_N x_N 		&\leq f_N  \\[-6pt]
    { {	{F_0 x_0}} }~~+~~  &{  F_1 x_1}~~+~~\cdots &~~+~~F_N x_N      	&= {b_{N+1}}  \\[-6pt]
    { {d_{N+1}}}\leq{{G_0 x_0}}~~+~~   &{  G_1 x_1}~~+~~\cdots &~~+~~G_N x_N			& \leq {f_{N+1}} \\[-6pt]
l_i \leq x_i &\leq u_i \quad \forall i=0,\dots,N.
\end{alignedat}
\end{equation}

The system matrix is split into sub-matrices $A_i \in \R^{m_{i_A}\times n_0}$, $B_i\in\R^{m_{i_A}\times n_i}$, $F_i \in \R^{m_{{N+1}_A}\times n_i}$, for the equality constraints and $C_i\in\R^{m_{i_C}\times n_0}$, $D_i\in\R^{m_{i_C}\times n_i}$, $G_i \in \R^{m_{{N+1}_C}\times n_i}$ for the inequality constraints, where $i \in \{0,\dots,N\}$. The matrices $F_i$, $G_i$ correspond to the linking constraints and the variables $x_0$ correspond to the linking variables. We call constraints associated with $\begin{bMatrix}{cc}A_i & B_i\end{bMatrix}$ or $\begin{bMatrix}{cc}C_i & D_i\end{bMatrix}$ \emph{local constraints}, and variables associated with $x_i$, $i\neq 0$ \emph{local variables}.

PIPS-IPM++ exploits this arrowhead structure to parallelize the linear algebra within its IPM and presolving. Within its IPM, PIPS-IPM++ uses MPI (Section~\ref{subsec:mpi}) to implement a Schur complement decomposition \cite{Rehfeldt2019_PIPSIPMpp, kempke2024massivelyparallelinteriorpointmethodarrowhead} to factorize the system matrix for each IPM iteration. To achieve this, for a given set of MPI processes (or simply processes), each diagonal block $i$ is assigned to exactly one process along with the block's data ($A_i$, $B_i$, $C_i$, $D_i$, $F_i$, $G_i$, $b_i$, $d_i$, $f_i$, $c_i$, $l_i$, $u_i$). Additionally, each process has access to the $0$-block data. This data-dependent distribution naturally limits the amount of MPI processes that can be employed to solve a problem to $N$, the number of diagonal blocks. Building on this distribution of data, each process performs a series of computations in parallel interleaved with serial MPI communication with the other processes. Should there be fewer processes available than blocks, multiple blocks are sequentially assigned to a single process. In the following presentation we assume that the number of processes and the number of blocks coincide.

This leads to a simple classification of presolve reductions. \emph{Local presolve reductions}, e.g., reductions considering only local constraints and variables, can be applied independently by each process and require no communication to be detected. Reductions on linking constraints or variables, \emph{global presolve reductions}, require communication among processes. An invariant enforced during our presolve is that no presolver may destroy the arrowhead structure, even if it could improve size or stability, to preserve compatibility with subsequent presolve rounds and PIPS-IPM++’s IPM. We also do not allow the creation of new linking constraints or linking variables, since increased coupling between blocks reduces block separability, increases MPI communication during presolve and solution, and thereby weakens the parallel scalability of the solver and presolver, which relies on large, weakly coupled diagonal blocks.

\subsection{Parallel postsolving}

Applying presolve to an LP produces the, often smaller, \emph{reduced problem}, which can subsequently be solved using any LP solver. This yields a solution to the reduced problem. To obtain a solution to the original problem, a procedure called \emph{postsolve} is applied. Depending on the solver and solution algorithm, one might obtain a primal solution satisfying Equation~(\ref{eq:pips_lp}), a dual solution satisfying Equation~(\ref{eq:pips_lp_dual}), or a primal-dual optimal solution satisfying both. Our postsolve is designed for a primal-dual IPM \cite{wright1997_primal_dual_ipm} and thus aims to recover a fully primal-dual feasible solution (including strict primal-dual complementary slackness).

Given an optimal solution to the reduced problem, postsolve reverts all reductions applied during presolve in reverse order. This is typically achieved using a stack. During presolve, whenever a reduction is applied, all data required to recover an optimal primal-dual solution during postsolve is pushed onto the stack. During postsolve, the stack is traversed in reverse order.

As mentioned previously, in our distributed parallel framework reductions may be either local or global, the latter affecting linking variables or constraints and requiring MPI communication. We maintain $N$ stacks, one per process. During presolve, each process pushes data from local and global reductions onto the stack. During postsolve, each process then traverses its own stack reverting local and global reductions. Local reductions can typically be reverted without communicating with other processes. Global reductions, such as removing a linking constraint or variable, typically involve MPI calls during their reversion. As in our framework MPI communication calls are blocking (Section~\ref{subsec:mpi}), global reductions act as synchronization points during both presolve and postsolve. In this way, postsolve is fully parallel where possible, and synchronized when necessary.

\subsection{Arrowhead-specific presolving techniques}

We implemented two presolvers specialized to arrowhead LPs: \textit{Permutation} and \textit{LinDependencies}, described below.

\subsubsection{Permutation Presolver}\label{sec:permutation}

The permutation presolver ensures that constraints and variables are in the most suitable position according to the AHLP in Equation~(\ref{eq:pips_blocklp}). For example, an equality linking constraint in $\begin{bMatrix}{cccc}F_0 & F_1 & \dots & F_N\end{bMatrix}$ should be non-empty in at least two $F_i$, $i>0$. Otherwise, it should be placed in $\begin{bMatrix}{cc}A_i & B_i\end{bMatrix}$ or $A_0$. Similarly, a linking variable should be non-empty in $A_i$ or $C_i$ for at least two $i\in\{1,\dots,N\}$. Otherwise, it belongs in the respective local block.

When the problem is read from file, constraints and variables may not be in their correct block position, and they may leave it during presolve as other constraints and variables are removed. To maintain correct constraint and variable placement, the permutation presolver is called each presolve round. We do not iterate the permutations to a fixed point: permutations that become possible only after other permutations, or after other presolvers removed constraints and variables, are applied in the next round. The permutation presolver is called once more at the end of presolve. The application of permutations also simplifies the implementation of other presolve reductions which can rely on the correct positioning of all constraints and variables. E.g., the SingletonConstr will look for singleton constraints only in the diagonal blocks, omitting all $A_i$, $F_i$, $C_i$, and $G_i$.

The constraint permutations that may be performed are displayed in Figure~\ref{fig:permutation_rows}. Local constraints, if empty in $B_i$ and $D_i$, are moved to $A_0$ and $C_0$ (Figure~\ref{fig:permutation_rows}, left). Linking constraints can become either local constraints (Figure~\ref{fig:permutation_rows}, middle), if nonzero only in $F_0$, $G_0$, and $F_i$ and $G_i$ for exactly one $i > 0$, or they are moved to the $A_0$ and $C_0$ (Figure~\ref{fig:permutation_rows}, right), if they are nonzero only in $F_0$ and $G_0$. Within one call, columns are permuted before constraints, since moving a constraint into $A_0$ or $C_0$ could prevent the permutation of a column to a local block. For each of the two steps, all candidates are first detected on the current problem and then moved at once. 

The possible variable permutations are shown in Figure~\ref{fig:permutation_cols}. A linking variable can become a local variable (Figure~\ref{fig:permutation_cols}, left), if it has nonzeros in $A_0$, $C_0$, and $A_i$, $C_i$ for exactly one $i > 0$. Local variables, which are empty in their respective $B_i$ and $D_i$ (Figure~\ref{fig:permutation_cols}, right) are permuted into $F_0$ and $G_0$. This last reduction is particularly important for PIPS-IPM++, since the variables that appear only in the linking part of the problem, $A_0$, $F_0$, $C_0$, and $G_0$, receive special treatment when forming the Schur complement. This is because such variables require no additional MPI communication when forming the Schur complement and appear sparse in the Schur complement matrix. We put a longer description of this mechanism into 
\ifzibreport
Appendix~\ref{appendix:0link}
\else
the online supplements, Appendix~A,
\fi
as this has not previously been documented in the literature.

\newcommand{\permsize}{0.36}

\newcommand{\rowpermbody}{
%%%% FIRST MOVE %%%%
\begin{tikzpicture}[scale=\permsize]
	% A blocks:
	\draw[gray, fill=gray!20] (0,5.5) rectangle +(1,1);
	\draw[gray, fill=gray!20] (0,3.3) rectangle +(1,1);
	\draw[gray, fill=gray!20] (0,1.1) rectangle +(1,1);
	\draw[gray, fill=gray!20] (0,0) rectangle +(1,1);
	% Blmat blocks:
	\draw[gray, fill=gray!20] (2.2,0) rectangle +(1,1);
	\draw[gray, fill=gray!20] (4.4,0) rectangle +(1,1);
	% B blocks:
	\draw[gray, fill=gray!20] (2.2,3.3) rectangle +(1,1);
	\draw[gray, fill=gray!20] (4.4,1.1) rectangle +(1,1);
	% dots between:
	%vertikal
	\path (0.5,2.3) -- node[auto=false]{\scalebox{1}[0.8]{$\vdots$}} (0.5,3.3);
	\path (0.5,4.4) -- node[auto=false]{\scalebox{1}[0.8]{$\vdots$}} (0.5,5.5);
	% horizontal
	\path (1.1,0.5) -- node[auto=false]{\scalebox{0.8}[1]{$\ldots$}} (2.2,0.5);
	\path (3.3,0.5) -- node[auto=false]{\scalebox{0.8}[1]{$\ldots$}} (4.4,0.5);	
	% schräg1
	\path (3.1,3.8) -- node[auto=false]{$\ddots$} (4.4,2.2);	
	\path (1.1,5.5) -- node[auto=false]{$\ddots$} (2.2,4.4);	
	% a_ik element and column k:
	\draw[black, fill=black] (0,3.1) rectangle +(1,0.2);
\end{tikzpicture}
\begin{tikzpicture}[scale=\permsize, baseline=-0.5ex]
	\draw[->, line width=1.5pt, color=black!70!black] 
	(0,2.7) -- (2,2.7);
\end{tikzpicture}
\begin{tikzpicture}[scale=\permsize]
	% A blocks:
	\draw[gray, fill=gray!20] (0,5.5) rectangle +(1,1);
	\draw[gray, fill=gray!20] (0,3.3) rectangle +(1,1);
	\draw[gray, fill=gray!20] (0,1.1) rectangle +(1,1);
	\draw[gray, fill=gray!20] (0,0) rectangle +(1,1);
	% Blmat blocks:
	\draw[gray, fill=gray!20] (2.2,0) rectangle +(1,1);
	\draw[gray, fill=gray!20] (4.4,0) rectangle +(1,1);
	% B blocks:
	\draw[gray, fill=gray!20] (2.2,3.3) rectangle +(1,1);
	\draw[gray, fill=gray!20] (4.4,1.1) rectangle +(1,1);
	% dots between:
	%vertikal
	\path (0.5,2.4) -- node[auto=false]{\scalebox{1}[0.8]{$\vdots$}} (0.5,3.3);
	\path (0.5,4.4) -- node[auto=false]{\scalebox{1}[0.8]{$\vdots$}} (0.5,5.5);
	% horizontal
	\path (1.1,0.5) -- node[auto=false]{\scalebox{0.8}[1]{$\ldots$}} (2.2,0.5);
	\path (3.3,0.5) -- node[auto=false]{\scalebox{0.8}[1]{$\ldots$}} (4.4,0.5);	
	% schräg
	\path (3.1,3.8) -- node[auto=false]{$\ddots$} (4.4,2.2);	
	\path (1.1,5.5) -- node[auto=false]{$\ddots$} (2.2,4.4);	
	% a_ik element and column k:
	\draw[black, fill=black] (0,5.3) rectangle +(1,0.2);
\end{tikzpicture}
\hspace{0.2cm}
\vrule width 1pt \hspace{3pt}
\hspace{0.2cm}
%%%% SECOND MOVE %%%%
\begin{tikzpicture}[scale=\permsize]
	% A blocks:
	\draw[gray, fill=gray!20] (0,5.5) rectangle +(1,1);
	\draw[gray, fill=gray!20] (0,3.3) rectangle +(1,1);
	\draw[gray, fill=gray!20] (0,1.1) rectangle +(1,1);
	\draw[gray, fill=gray!20] (0,0) rectangle +(1,1);
	% Blmat blocks:
	\draw[gray, fill=gray!20] (2.2,0) rectangle +(1,1);
	\draw[gray, fill=gray!20] (4.4,0) rectangle +(1,1);
	% B blocks:
	\draw[gray, fill=gray!20] (2.2,3.3) rectangle +(1,1);
	\draw[gray, fill=gray!20] (4.4,1.1) rectangle +(1,1);
	% dots between:
	%vertikal
	\path (0.5,2.4) -- node[auto=false]{\scalebox{1}[0.8]{$\vdots$}} (0.5,3.3);
	\path (0.5,4.4) -- node[auto=false]{\scalebox{1}[0.8]{$\vdots$}} (0.5,5.5);
	% horizontal
	\path (1.1,0.5) -- node[auto=false]{\scalebox{0.8}[1]{$\ldots$}} (2.2,0.5);
	\path (3.3,0.5) -- node[auto=false]{\scalebox{0.8}[1]{$\ldots$}} (4.4,0.5);	
	% schräg
	\path (3.1,3.8) -- node[auto=false]{$\ddots$} (4.4,2.2);	
	\path (1.1,5.5) -- node[auto=false]{$\ddots$} (2.2,4.4);	
	% a_ik element and column k:
	\draw[black, fill=black] (0,-0.2) rectangle +(1,0.2);
	\draw[black, fill=black] (2.2,-0.2) rectangle +(1,0.2);
\end{tikzpicture}
\begin{tikzpicture}[scale=\permsize, baseline=-0.5ex]
	\draw[->, line width=1.5pt, color=black!70!black] 
	(0,2.7) -- (2,2.7);
\end{tikzpicture}
\begin{tikzpicture}[scale=\permsize]
	% A blocks:
	\draw[gray, fill=gray!20] (0,5.5) rectangle +(1,1);
	\draw[gray, fill=gray!20] (0,3.3) rectangle +(1,1);
	\draw[gray, fill=gray!20] (0,1.1) rectangle +(1,1);
	\draw[gray, fill=gray!20] (0,0) rectangle +(1,1);
	% Blmat blocks:
	\draw[gray, fill=gray!20] (2.2,0) rectangle +(1,1);
	\draw[gray, fill=gray!20] (4.4,0) rectangle +(1,1);
	% B blocks:
	\draw[gray, fill=gray!20] (2.2,3.3) rectangle +(1,1);
	\draw[gray, fill=gray!20] (4.4,1.1) rectangle +(1,1);
	% dots between:
	%vertikal
	\path (0.5,2.3) -- node[auto=false]{\scalebox{1}[0.8]{$\vdots$}} (0.5,3.3);
	\path (0.5,4.4) -- node[auto=false]{\scalebox{1}[0.8]{$\vdots$}} (0.5,5.5);
	% horizontal
	\path (1.1,0.5) -- node[auto=false]{\scalebox{0.8}[1]{$\ldots$}} (2.2,0.5);
	\path (3.3,0.5) -- node[auto=false]{\scalebox{0.8}[1]{$\ldots$}} (4.4,0.5);	
	% schräg
	\path (3.3,3.6) -- node[auto=false]{$\ddots$} (4.4,2.2);	
	\path (1.1,5.5) -- node[auto=false]{$\ddots$} (2.2,4.4);	
	% a_ik element and column k:
	\draw[black, fill=black] (0,3.1) rectangle +(1,0.2);
	\draw[black, fill=black] (2.2,3.1) rectangle +(1,0.2);
\end{tikzpicture}
\hspace{0.2cm}
\vrule width 1pt \hspace{3pt}
\hspace{0.2cm}
%%%% THIRD MOVE %%%%
\begin{tikzpicture}[scale=\permsize]
	% A blocks:
	\draw[gray, fill=gray!20] (0,5.5) rectangle +(1,1);
	\draw[gray, fill=gray!20] (0,3.3) rectangle +(1,1);
	\draw[gray, fill=gray!20] (0,1.1) rectangle +(1,1);
	\draw[gray, fill=gray!20] (0,0) rectangle +(1,1);
	% Blmat blocks:
	\draw[gray, fill=gray!20] (2.2,0) rectangle +(1,1);
	\draw[gray, fill=gray!20] (4.4,0) rectangle +(1,1);
	% B blocks:
	\draw[gray, fill=gray!20] (2.2,3.3) rectangle +(1,1);
	\draw[gray, fill=gray!20] (4.4,1.1) rectangle +(1,1);
	% dots between:
	%vertikal
	\path (0.5,2.4) -- node[auto=false]{\scalebox{1}[0.8]{$\vdots$}} (0.5,3.3);
	\path (0.5,4.4) -- node[auto=false]{\scalebox{1}[0.8]{$\vdots$}} (0.5,5.5);
	% horizontal
	\path (1.1,0.5) -- node[auto=false]{\scalebox{0.8}[1]{$\ldots$}} (2.2,0.5);
	\path (3.3,0.5) -- node[auto=false]{\scalebox{0.8}[1]{$\ldots$}} (4.4,0.5);	
	% schräg
	\path (3.1,3.8) -- node[auto=false]{$\ddots$} (4.4,2.2);	
	\path (1.1,5.5) -- node[auto=false]{$\ddots$} (2.2,4.4);	
	% a_ik element and column k:
	\draw[black, fill=black] (0,-0.2) rectangle +(1,0.2);
\end{tikzpicture}
\begin{tikzpicture}[scale=\permsize, baseline=-0.5ex]
	\draw[->, line width=1.5pt, color=black!70!black] 
	(0,2.7) -- (2,2.7);
\end{tikzpicture}
\begin{tikzpicture}[scale=\permsize]
	% A blocks:
	\draw[gray, fill=gray!20] (0,5.5) rectangle +(1,1);
	\draw[gray, fill=gray!20] (0,3.3) rectangle +(1,1);
	\draw[gray, fill=gray!20] (0,1.1) rectangle +(1,1);
	\draw[gray, fill=gray!20] (0,0) rectangle +(1,1);
	% Blmat blocks:
	\draw[gray, fill=gray!20] (2.2,0) rectangle +(1,1);
	\draw[gray, fill=gray!20] (4.4,0) rectangle +(1,1);
	% B blocks:
	\draw[gray, fill=gray!20] (2.2,3.3) rectangle +(1,1);
	\draw[gray, fill=gray!20] (4.4,1.1) rectangle +(1,1);
	% dots between:
	%vertikal
	\path (0.5,2.4) -- node[auto=false]{\scalebox{1}[0.8]{$\vdots$}} (0.5,3.3);
	\path (0.5,4.4) -- node[auto=false]{\scalebox{1}[0.8]{$\vdots$}} (0.5,5.5);
	% horizontal
	\path (1.1,0.5) -- node[auto=false]{\scalebox{0.8}[1]{$\ldots$}} (2.2,0.5);
	\path (3.3,0.5) -- node[auto=false]{\scalebox{0.8}[1]{$\ldots$}} (4.4,0.5);	
	% schräg
	\path (3.1,3.8) -- node[auto=false]{$\ddots$} (4.4,2.2);	
	\path (1.1,5.5) -- node[auto=false]{$\ddots$} (2.2,4.4);	
	% a_ik element and column k:
	\draw[black, fill=black] (0,5.3) rectangle +(1,0.2);
\end{tikzpicture}
}

\ifzibreport
\begin{figure}[t]
\caption{Constraint permutations executed by the permutation presolver.}\label{fig:permutation_rows}
\adjustbox{max width=\textwidth}{%
\rowpermbody
}
\end{figure}
\else
\begin{figure}[t]
\FIGURE
{\rowpermbody}
{Constraint permutations executed by the permutation presolver.\label{fig:permutation_rows}}
{}
\end{figure}
\fi

\newcommand{\colpermbody}{
\centering
%%%% FIRST MOVE %%%%
\begin{tikzpicture}[scale=\permsize]
	% A blocks:
	\draw[gray, fill=gray!20] (0,5.5) rectangle +(1,1);
	\draw[gray, fill=gray!20] (0,3.3) rectangle +(1,1);
	\draw[gray, fill=gray!20] (0,1.1) rectangle +(1,1);
	\draw[gray, fill=gray!20] (0,0) rectangle +(1,1);
	% Blmat blocks:
	\draw[gray, fill=gray!20] (2.2,0) rectangle +(1,1);
	\draw[gray, fill=gray!20] (4.4,0) rectangle +(1,1);
	% B blocks:
	\draw[gray, fill=gray!20] (2.2,3.3) rectangle +(1,1);
	\draw[gray, fill=gray!20] (4.4,1.1) rectangle +(1,1);
	% dots between:
	%vertikal
	\path (0.5,2.4) -- node[auto=false]{\scalebox{1}[0.8]{$\vdots$}} (0.5,3.3);
	\path (0.5,4.4) -- node[auto=false]{\scalebox{1}[0.8]{$\vdots$}} (0.5,5.5);
	% horizontal
	\path (1.1,0.5) -- node[auto=false]{\scalebox{0.8}[1]{$\ldots$}} (2.2,0.5);
	\path (3.3,0.5) -- node[auto=false]{\scalebox{0.8}[1]{$\ldots$}} (4.4,0.5);	
	% schräg1
	\path (3.1,3.8) -- node[auto=false]{$\ddots$} (4.4,2.2);	
	\path (1.1,5.5) -- node[auto=false]{$\ddots$} (2.2,4.4);	
	% a_ik element and column k:
	\draw[black, fill=black] (1,3.3) rectangle +(0.2,1);
	\draw[black, fill=black] (1,0) rectangle +(0.2,1);
\end{tikzpicture}
\begin{tikzpicture}[scale=0.36, baseline=-0.5ex]
	\draw[->, line width=1.5pt, color=black!70!black] 
	(0,2.7) -- (2,2.7);
\end{tikzpicture}
\begin{tikzpicture}[scale=\permsize]
	% A blocks:
	\draw[gray, fill=gray!20] (0,5.5) rectangle +(1,1);
	\draw[gray, fill=gray!20] (0,3.3) rectangle +(1,1);
	\draw[gray, fill=gray!20] (0,1.1) rectangle +(1,1);
	\draw[gray, fill=gray!20] (0,0) rectangle +(1,1);
	% Blmat blocks:
	\draw[gray, fill=gray!20] (2.2,0) rectangle +(1,1);
	\draw[gray, fill=gray!20] (4.4,0) rectangle +(1,1);
	% B blocks:
	\draw[gray, fill=gray!20] (2.2,3.3) rectangle +(1,1);
	\draw[gray, fill=gray!20] (4.4,1.1) rectangle +(1,1);
	% dots between:
	%vertikal
	\path (0.5,2.4) -- node[auto=false]{\scalebox{1}[0.8]{$\vdots$}} (0.5,3.3);
	\path (0.5,4.4) -- node[auto=false]{\scalebox{1}[0.8]{$\vdots$}} (0.5,5.5);
	% horizontal
	\path (1.1,0.5) -- node[auto=false]{\scalebox{0.8}[1]{$\ldots$}} (2.2,0.5);
	\path (3.3,0.5) -- node[auto=false]{\scalebox{0.8}[1]{$\ldots$}} (4.4,0.5);	
	% schräg
	\path (3.1,3.8) -- node[auto=false]{$\ddots$} (4.4,2.2);	
	\path (1.1,5.5) -- node[auto=false]{$\ddots$} (2.2,4.4);	
	% a_ik element and column k:
	\draw[black, fill=black] (3.2,3.3) rectangle +(0.2,1);
	\draw[black, fill=black] (3.2,0) rectangle +(0.2,1);
\end{tikzpicture}
\hspace{0.2cm}
\vrule width 1pt \hspace{3pt}
\hspace{0.2cm}
%%%% SECOND MOVE %%%%
\begin{tikzpicture}[scale=\permsize]
	% A blocks:
	\draw[gray, fill=gray!20] (0,5.5) rectangle +(1,1);
	\draw[gray, fill=gray!20] (0,3.3) rectangle +(1,1);
	\draw[gray, fill=gray!20] (0,1.1) rectangle +(1,1);
	\draw[gray, fill=gray!20] (0,0) rectangle +(1,1);
	% Blmat blocks:
	\draw[gray, fill=gray!20] (2.2,0) rectangle +(1,1);
	\draw[gray, fill=gray!20] (4.4,0) rectangle +(1,1);
	% B blocks:
	\draw[gray, fill=gray!20] (2.2,3.3) rectangle +(1,1);
	\draw[gray, fill=gray!20] (4.4,1.1) rectangle +(1,1);
	% dots between:
	%vertikal
	\path (0.5,2.4) -- node[auto=false]{\scalebox{1}[0.8]{$\vdots$}} (0.5,3.3);
	\path (0.5,4.4) -- node[auto=false]{\scalebox{1}[0.8]{$\vdots$}} (0.5,5.5);
	% horizontal
	\path (1.1,0.5) -- node[auto=false]{\scalebox{0.8}[1]{$\ldots$}} (2.2,0.5);
	\path (3.3,0.5) -- node[auto=false]{\scalebox{0.8}[1]{$\ldots$}} (4.4,0.5);	
	% schräg
	\path (3.1,3.8) -- node[auto=false]{$\ddots$} (4.4,2.2);	
	\path (1.1,5.5) -- node[auto=false]{$\ddots$} (2.2,4.4);	
	% a_ik element and column k:
	\draw[black, fill=black] (3.2,0) rectangle +(0.2,1);
\end{tikzpicture}
\begin{tikzpicture}[scale=\permsize, baseline=-0.5ex]
	\draw[->, line width=1.5pt, color=black!70!black] 
	(0,2.7) -- (2,2.7);
\end{tikzpicture}
\begin{tikzpicture}[scale=\permsize]
	% A blocks:
	\draw[gray, fill=gray!20] (0,5.5) rectangle +(1,1);
	\draw[gray, fill=gray!20] (0,3.3) rectangle +(1,1);
	\draw[gray, fill=gray!20] (0,1.1) rectangle +(1,1);
	\draw[gray, fill=gray!20] (0,0) rectangle +(1,1);
	% Blmat blocks:
	\draw[gray, fill=gray!20] (2.2,0) rectangle +(1,1);
	\draw[gray, fill=gray!20] (4.4,0) rectangle +(1,1);
	% B blocks:
	\draw[gray, fill=gray!20] (2.2,3.3) rectangle +(1,1);
	\draw[gray, fill=gray!20] (4.4,1.1) rectangle +(1,1);
	% dots between:
	%vertikal
	\path (0.5,2.3) -- node[auto=false]{\scalebox{1}[0.8]{$\vdots$}} (0.5,3.3);
	\path (0.5,4.4) -- node[auto=false]{\scalebox{1}[0.8]{$\vdots$}} (0.5,5.5);
	% horizontal
	\path (1.1,0.5) -- node[auto=false]{\scalebox{0.8}[1]{$\ldots$}} (2.2,0.5);
	\path (3.3,0.5) -- node[auto=false]{\scalebox{0.8}[1]{$\ldots$}} (4.4,0.5);	
	% schräg
	\path (3.3,3.6) -- node[auto=false]{$\ddots$} (4.4,2.2);	
	\path (1.1,5.5) -- node[auto=false]{$\ddots$} (2.2,4.4);	
	% a_ik element and column k:
	\draw[black, fill=black] (1,0) rectangle +(0.2,1);
\end{tikzpicture}
}

\ifzibreport
\begin{figure}[t]
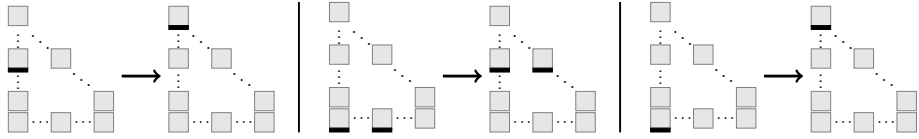

\caption{Variable permutations executed by the permutation presolver.}\label{fig:permutation_cols}
\colpermbody
\end{figure}
\else
\begin{figure}[t]
\FIGURE
{\colpermbody}
{Variable permutations executed by the permutation presolver.\label{fig:permutation_cols}}
{}
\end{figure}
\fi

\subsubsection{Linear Dependencies Presolver}\label{sec:lindep}

The LinDependencies presolver detects linearly dependent constraints in the equality matrix, ensuring full rank of the KKT system used in the IPM, which is critical for stability. Applying a full distributed Gaussian elimination for linking constraints however, would be prohibitively expensive. Instead, we detect dependencies only in the equality matrix excluding linking constraints. As such, it suffices to check that all $B_i$, $i=1,\dots,N$, and $A_0$ are of full rank to establish that the sub-matrix
\begin{equation*}
\begin{bMatrix}{cccc}
    A_0    &       &        &  \\
    A_1    & B_1   &        & \\
    \vdots &       & \ddots & \\
    A_N    &       &        & B_N \\
\end{bMatrix}
\end{equation*}
is of full rank. These checks can be performed independently. The procedure thus guarantees full row rank of the equality constraints without the linking constraints. Linear dependencies that involve linking constraints are not detected. We do not consider this a limitation in practice: the linking constraints do not enter the factorized matrices directly but only in aggregated form through the Schur complement, and detecting such dependencies would require either a distributed elimination or gathering all linking constraints on a single process, both of which require extensive communication. In our experiments, any remaining dependencies are left to the linear solver factorizing the Schur complement.

Each process $i$ constructs its local matrix $\begin{bMatrix}{cc}B_i & I\end{bMatrix}$, with $I$ a unit matrix of size $m_{i_A}\times m_{i_A}$, and applies a Gaussian elimination algorithm similar to MA50, the factorization kernel of MA48~\cite{Duff1993_MA48}. Pivoting is restricted to the first $n_i$ matrix columns. After elimination, the matrix is partitioned as
\begin{equation*}
\begin{bMatrix}{cc}
    \hat{B}_i & K_{1} \\
    0         & K_{2} \\
\end{bMatrix}
\end{equation*}
where $\hat{B}_i \in \R^{\hat{m}_{i_A}\times n_i}$, $K_1 \in \R^{\hat{m}_{i_A}\times m_{i_A}}$, and $K_2 \in \R^{(m_{i_A}-\hat{m}_{i_A})\times m_{i_A}}$. Each row in $K_1$ and $K_2$ represents the linear combination applied to the corresponding constraint of $\begin{bMatrix}{cc}A_i & B_i\end{bMatrix}$:
\begin{equation*}
\begin{bMatrix}{c} K_1 \\ K_2 \\\end{bMatrix} \begin{bMatrix}{cc}A_i & B_i \\\end{bMatrix} = P\begin{bMatrix}{cc}\hat{A}_{i1} & \hat{B}_i \\ \hat{A}_{i2} & 0 \\\end{bMatrix},
\end{equation*}
where $P \in \N_0^{m_{i_A}\times m_{i_A}}$ is a permutation matrix corresponding to the row/constraint permutation applied during the elimination. Constraints where $\hat{A}_{i2}$ is empty are fully linearly dependent. Depending on the constraint right-hand sides they either prove infeasibility or are removed from the problem. Constraints where $\hat{A}_{i2}$ is non-empty are moved to $A_0$. Finally, all processes jointly apply Gaussian elimination to $A_0$ to detect further dependencies. The linear combinations used during elimination are stored on the stack for proper reversal during postsolve.
\section{Algorithms and software design}\label{sec:algos}

The communication of data during an individual presolver is critical for the efficiency of the presolve procedure and data is communicated both during the presolve and postsolve procedure.
This is particularly important when performing pre- and postsolve procedures on linking constraints and variables.
This section will first describe the communication of data in pre- and postsolve.
A description of the parallel constraints presolver is then provided to illustrate the required communication among processes and the differences between local and global reductions.

\subsection{Data communication during presolve}

Communication during our presolve happens in two ways: during the application of a presolver, by reducing, e.g., communicating a hash as in the parallel constraint detection (Section~\ref{subsec:parallel_rows}), or by periodically communicating problem characteristics as explained in the following. To apply our presolve efficiently without recomputing this data frequently, our implementation keeps track of the nonzeros of constraints and variables in the problem as well as the constraint activities, where the minimum (actmin) and maximum activity (actmax) of an (equality) constraint $i$ is given as
\ifzibreport
\begin{equation*}
	\text{actmin}_i = \min_{l \leq x \leq u} A_{i.}\T x = \sum_{a_{ij} > 0} l_j a_{ij} + \sum_{a_{ij} < 0} u_j a_{ij},
\end{equation*}
\begin{equation*}
	\text{actmax}_i = \max_{l \leq x \leq u} A_{i.}\T x = \sum_{a_{ij} < 0} l_j a_{ij} + \sum_{a_{ij} > 0} u_j a_{ij}.
\end{equation*}
\else
\begin{equation*}
	\text{actmin}_i = \min_{l \leq x \leq u} A_{i.}\T x = \sum_{a_{ij} > 0} l_j a_{ij} + \sum_{a_{ij} < 0} u_j a_{ij},\, \text{actmax}_i = \max_{l \leq x \leq u} A_{i.}\T x = \sum_{a_{ij} < 0} l_j a_{ij} + \sum_{a_{ij} > 0} u_j a_{ij}.
\end{equation*}
\fi
Actmin and actmax are defined equivalently for inequality constraints. Should any of the variables used in either of the activity formulas be unbounded, the respective activity is set to $\pm \inf$. Whenever a change in the problem occurs, these quantities get updated.

For linking variables and constraints, we cannot directly update the quantities, as changes might result from local reductions, known only to a single process. Rather, such changes are buffered and communicated periodically among all processes, typically when a presolver finishes its work. 

\subsection{Data communication during postsolve}\label{sec:comm_postsolve}

During postsolve, local reversions of presolve transformations must be consistent with globally shared quantities, such as dual variables and activities associated with linking constraints and variables. This requires MPI communication to ensure that all processes maintain a consistent view of these global quantities after local postsolve operations have been applied.

Analogous to presolve, postsolve communication occurs in two ways. First, communication is required during a postsolve operation whenever a quantity that depends on distributed data must be evaluated, for example when computing the activity of a linking constraint. Second, communication is required to synchronize buffered changes that originate from local postsolve operations and affect globally shared data.

The latter case arises, for instance, when a dual value is shifted between a constraint and a variable to revert a bound tightening identified during presolve. If the constraint is local but affects the bound of a linking variable, the corresponding dual of the linking variable is modified by a local reduction. Reverting this bound tightening on a single process leaves the dual values of the linking variable outdated on all other processes. To address this, we employ a buffering mechanism analogous to the one used for tracking global data during presolve: instead of immediately applying changes to globally required data (e.g., duals of linking constraints or variables), local modifications are buffered.

During presolve, \emph{synchronization events} are recorded on the stack of each process whenever such globally relevant updates are deferred. In postsolve, after the corresponding local reversions have been applied, the stored synchronization event is executed and all buffered local changes are communicated to restore global consistency.

\subsection{Parallel constraint detection}\label{subsec:parallel_rows}

The parallel constraint detection in PIPS-IPM++ works similarly to the one described in \cite{AchterbergBixbyGuetal2016}. We apply a two level hashing algorithm to first identify constraints with the same support and later, within each hash bucket of constraints with the same support, identify constraints with the same coefficients. As is usually the case in our framework, there are two different implementations, one for local constraints and one for global/linking constraints. Local parallel constraints are detected independently within each block. For linking constraints, each process computes the support and coefficient hashes with respect to its local blocks $F_i$ and $G_i$. The processes combine the support hashes by an \texttt{Allreduce} with a custom MPI operator, and the information whether a pair of constraints is parallel locally by an \texttt{Allreduce} with the logical AND operator, after which each process applies the same reductions. We do not try to detect parallel local and global constraints, but rely on the permutation presolver to first put these constraints into their correct positions (see Section~\ref{sec:permutation}). Managing this communication requires careful design and implementation to balance the trade-off between communication overhead and reduction strength, occasionally weakening presolve to limit overhead. We give pseudocode and a detailed description of both implementations in
\ifzibreport
Appendix~\ref{appendix:parallel_rows}.
\else
the online supplement, Appendix~C.
\fi

\section{Computational Experiments} \label{sec:results}

We conducted four experiments to evaluate our presolve implementation. In Section~\ref{sec:results_scaling}, we assess the scalability of our presolve on eight large instances. In Section~\ref{sec:results_efficiency} we compare our implementation with respect to runtime and number of reductions with the presolve implementations of Gurobi and PaPILO. In Section~\ref{sec:results_impact_ipm} we evaluate the impact of our presolving on running the IPM implemented in PIPS-IPM++. Finally, in Section~\ref{sec:results_external} we compare our presolve with presolving the model externally with Gurobi or PaPILO before passing it to PIPS-IPM++. For the evaluation of our experiments, we generally use the shifted geometric mean (sgm) to lessen the impact of extreme outliers, using a shift of one second and a shift of one percent for the respective aggregated results.

\paragraph{Testset} PIPS-IPM++, and thus our presolve implementation, can currently only be called via a GAMS (\url{https://www.gams.com}) interface. This interface requires models to be annotated to assign variables and constraints to an arrowhead form (see Equation~(\ref{eq:pips_blocklp})). Our model library contains 82 instances of varying sizes, most of which come from the energy system context. We note that not all models in our testset expose a block structure that is favorable towards PIPS-IPM++'s Schur complement approach. As such, PIPS-IPM++ does not always outperform commercial solvers such as Gurobi on our testset. We have also used instances that PIPS-IPM++ cannot solve at all due to an overly large Schur complement. Still, we included these models to get a broad assessment of our presolve implementation. Whenever we run Gurobi and PaPILO on instances from our testset, we first converted the annotated GAMS files to MPS files with GAMS convert, or with PIPS-IPM++ for the instances that are available only as block files. We display sizes, number of blocks, and linking part of all instances in
\ifzibreport
Appendix~\ref{appendix:models}, Table~\ref{tab:sizes}.
\else
the online supplements, Appendix~B, Table~S1.
\fi

\paragraph{Hardware} All experiments were conducted on the terrabyte supercomputer of the Leibniz Supercomputing Centre. Each node is equipped with 1024 GB RAM and two Intel Xeon Platinum 8380 CPUs with 40 cores each at 2.3 GHz (80 cores and 160 hardware threads per node). All jobs always allocated full nodes, so no concurrent execution was permitted.

\paragraph{Software} We used PaPILO 3.0.1, Gurobi 13.0.0 and PIPS-IPM++ (git hash \texttt{d34cee38}).

\paragraph{Timing} All reported PIPS-IPM++ times are wall-clock times from the start of the program (after MPI initialization) until postsolve has finished, including reading the problem, presolve, the IPM, and postsolve. Each process reads only the data of its own blocks and of the 0-block from the block-wise GDX files, so there is no separate distribution step. Since the 0-block is available on every process, postsolve does not need to move solution values of linking variables and constraints. Model generation time is not included. On our testset, reading takes 0.98\,s and postsolve 0.16\,s (sgm), and reading, presolve, and postsolve take a median of 0.1\,\%, 0.4\,\%, and 0.02\,\% of the total time. For Gurobi and PaPILO we report the time of the presolve alone, excluding reading the MPS file.

\paragraph{Reproducibility} All results can be reproduced with PIPS-IPM++ git hash \texttt{d34cee38}. Some runs were made with earlier versions that differ from it only in two bug fixes which do not affect these runs. Apart from an IPM time limit of one hour and Panua-Pardiso as linear solver, PIPS-IPM++ ran with its default settings. All runs of PIPS-IPM++ are deterministic for a fixed number of processes: we enforce a fixed order of the reduction operations within MPI collectives with Intel MPI's conditional numerical reproducibility (\texttt{I\_MPI\_CBWR=2}). Gurobi's presolve was called through its Python interface with a time limit of one hour, PaPILO's presolve with its default settings and a job limit of 1.5 hours, both on MPS files exported from the annotated models.

\paragraph{Data availability} 69 of the 82 instances are available at \url{https://doi.org/10.5281/zenodo.23269066} as annotated GDX files, which contain the block assignment of constraints and variables and are read by PIPS-IPM++, and as MPS files. The remaining instances were provided by project partners and cannot be published. The record contains the scripts used to run and evaluate the experiments and the log files of runs on these instances. Instance names follow Table~S1 in the online supplement.

\paragraph{Scope and limitations} Our testset consists almost entirely of energy system models, annotated by the modelers (for the SIMPLE instances through the model generator). As Section~\ref{sec:results_external} shows, the quality of this annotation is crucial, and an automatic annotation cannot yet replace it. Of the presolvers in Table~\ref{tab:presolvers}, only Permutation and LinDependencies are specific to the arrowhead structure. We expect our results to carry over to other AHLPs with a good block structure, but we have not tested this. In our testset, the blocks are almost equally large (on the 62 instances with a complete annotation, the largest block has at most 8\% more variables than the average block). PIPS-IPM++ assigns consecutive blocks to the processes such that the number of blocks per process differs by at most one, it does not balance their size. With uneven blocks, the process holding the largest blocks determines the time of the local reductions. Instances with large linking parts gain little from more than one compute node (Section~\ref{sec:results_efficiency}).

\subsection{Scaling} \label{sec:results_scaling}

First, we demonstrate the scaling of our presolve routines on eight large instances: all instances of our testset with at least 200 million nonzeros and at least 96 blocks, except for variants of nearly the same size (elmod\_cwe14 and elmod\_cwe16 of elmod\_cwe15, yssp\_exp\_lim\_488r\_250b of yssp\_exp\_488r\_250b), and eustem (460 million nonzeros, 70 blocks). The instance yssp\_exp\_488r\_1460b is yssp\_exp\_488r\_250b annotated with 1460 instead of 250 blocks. It and eustem were chosen as they have comparatively large linking parts (\ifzibreport Table~\ref{tab:sizes}\else online supplement, Table~S1\fi). We ran Gurobi and PaPILO with up to 128 threads, the maximum number that, allowing for hyperthreading, could still be run on a single compute node, and PIPS-IPM++ with one core per process and at most one process per block. Table~\ref{tab:scaling_details} gives size and block count of these instances and the nonzeros remaining after presolve (Section~\ref{sec:results_efficiency} analyzes the reductions on the whole testset). Gurobi removes 3 to 15 percentage points more nonzeros than our presolve, except on eustem, where our presolve removes slightly more than both Gurobi and PaPILO. PaPILO removes the most on elmod\_cwe15, did not finish presolving the three yssp\_exp instances within 1.5 hours, and aborted with an error on elmod\_set5.

\newcommand{\scalingrows}{%
\toprule
\multicolumn{5}{c}{} & \multicolumn{3}{c}{\% nonzeros remaining} \\
\cmidrule(lr){6-8} 
Instance  & \# blocks & \# nonzeros & \# constraints & \# variables & PIPS-IPM++ & Gurobi & PaPILO \\
\midrule
tb\_400\_10\_shift    & 1024 & 205\,569\,328 &  51\,604\,093 &  59\,795\,108 & 81.23 & 78.08 & 81.25 \\
%simple/TO_0974_TBSIZE_400_NBREGIONS_10_RESOLUTION_5_240_NBSHIFTS2_SHIFTSTEPSIZE_2/block_noVEnames
simple\_tudpaper1    &  438 & 207\,679\,112 &  52\,193\,591 &  57\,444\,984 & 91.29 & 87.35 & 91.36 \\
%beam-me-mext/simple_instances/SIMPLE_tudpaper1/allblocksPips
elmod\_cwe15     &  438 & 271\,875\,064 &  98\,646\,274 &  85\,646\,554 & 57.33 & 44.69 & 39.30 \\
%beam-me-mext/TUD/Set4_CWE15/ELMOD_438_20_noVEnames
elmod\_set5     &  438 & 711\,769\,260 & 254\,304\,960 & 224\,677\,685 & 58.40 & 45.64 & Error \\
%beam-me-mext/TUD/Set5/ELMOD_438_20_noVEnames
yssp\_exp\_dcopf\_488r\_96b &   96 & 316\,863\,066 &  96\,851\,394 & 110\,650\,876 & 54.07 & 43.47 & Timeout \\
%beam-me-mext/yssp/YSSP_exp_DCopf_TaW_96b/blocks_YSSP_exp_DCopf_TaW_488r_1h_t_96b_novenames
yssp\_exp\_488r\_250b &  250 & 247\,383\,863 &  73\,253\,433 &  87\,054\,963 & 69.15 & 53.82 & Timeout \\
%beam-me-mext/yssp/exp_488r_t_250b/YSSP_exp_488r_t_250b_novenames
yssp\_exp\_488r\_1460b & 1460 & 247\,383\,863 &  73\,253\,433 &  87\,054\,963 & 69.15 & 53.82 & Timeout \\
%beam-me-mext/yssp/exp_488r_t_1460b/YSSP_exp_488r_t_1460b_novenames
eustem       &   70 & 460\,037\,503 &  28\,635\,704 &  19\,195\,771 & 94.76 & 95.30 & 96.36 \\
%beam-me-mext/PSI/Set3/eustem_672_11_8jac
\bottomrule
}
\ifzibreport
\begin{table}[ht!]
\centering
\caption{Instances and results for scaling experiment.}
\label{tab:scaling_details}
\adjustbox{max width=\textwidth}{%
\begin{tabular}{l|rrrrrrr}
\scalingrows
\end{tabular}
}
\end{table}
\else
\begin{table}[ht!]
\TABLE
{Instances and results for scaling experiment.\label{tab:scaling_details}}
{\begin{tabular}{l|rrrrrrr}
\scalingrows
\end{tabular}}
{}
\end{table}
\fi
Figure~\ref{fig:presolve_scaling} plots the presolve time of PIPS-IPM++, Gurobi, and PaPILO against the number of MPI processes (PIPS-IPM++) or threads (Gurobi and PaPILO), together with the optimal linear speed-up as a reference.
Neither Gurobi nor PaPILO achieve a significant speed-up from additional threads on the instances given. Gurobi's presolve is essentially sequential on our instances (sgm 26.7\,s with one and 26.6\,s with 64 threads over the instances of Section~\ref{sec:results_efficiency}). PaPILO runs its presolvers in parallel in shared memory, but gains only a factor of 1.25 from one to 64 threads (sgm). On the other hand, PIPS-IPM++ performs worse than Gurobi, and except on elmod\_cwe15 and eustem also than PaPILO, when using only one thread, a fact that we attribute to the additional overhead required by the parallel presolve implementation as well as the AHLP treatment within our software. For the given instances, the breakeven point of PIPS-IPM++'s presolve time and Gurobi and PaPILO's lies mostly at two to four processes, and at 16 processes for yssp\_exp\_488r\_1460b. On the six instances with small linking parts, our presolve is 30 to 50 times faster with 64 processes than with one process, and 36 to 172 times faster with the maximum number of processes (a sequential fraction of 0.4 to 1.8\% at 64 processes in terms of Amdahl's law~\cite{Amdahl1967}). On eustem and yssp\_exp\_488r\_1460b, the speed-up at 64 processes is only 23 and 7. On the latter, our presolve is fastest with 16 processes and slower than Gurobi from 64 processes on. The amount of sequential processing within PIPS-IPM++ grows with the number of linking variables and constraints, whereas a fully block-diagonal matrix would require no communication and should scale linearly. Moreover, the cost per MPI operation increases with the number of processes, most strongly when going from one to multiple compute nodes (from 64 to 128 and more processes).   
\ifzibreport
\begin{figure}
	\captionsetup[sub]{font=small} % Only affects sub-captions
	\begin{subfigure}{.5\textwidth}
		\caption{tb\_400\_10\_shift}
		\begin{adjustbox}{width=0.95\textwidth}
			\begin{tikzpicture}
				\begin{loglogaxis}
					[
					xlabel={\MPI processes/threads},
					ylabel={Presolve time [s]},
					log ticks with fixed point,
					scale only axis,
					xmin=1, xmax=1024,
					ymin=0.1, ymax=1024,
					legend pos=south west,
					grid style={line width=1pt, draw=black, opacity=0.1},
					xmajorgrids={true},
					ymajorgrids={true},
					log basis x=2
					]
					\addplot %simple/TO_0974_TBSIZE_400_NBREGIONS_10_RESOLUTION_5_240_NBSHIFTS2_SHIFTSTEPSIZE_2/block_noVEnames
					[
					color=cp6,
					mark=square*,
					thick
					]
					coordinates {
						(1,413.932)
						(2,204.421)
						(4,99.801)
						(8,67.963)
						(16,34.595)
						(32,13.176)
						(64,11.075)
						(128,12.200)
						(256,5.913)
						(512,4.456)
						(1024,3.576)
					};
					\addplot %gurobi
					[
					color=cp5,
					mark=triangle*,
					thick
					]
					coordinates {
						(1,206.47)
						(2,207.03)
						(4,204.87)
						(8,207.46)
						(16,204.85)
						(32,207.41)
						(64,205.28)
						(128,205.76)
					};
					\addplot % papilo
					[
					color=cp2,
					mark=diamond*,
					thick
					]
					coordinates {
						(1,225.416)
						(2,176.424)
						(4,146.075)
						(8,171.193)
						(16,168.503)
						(32,170.547)
						(64,142.928)
						(128,183.375)
					};
					\addplot % ideal speedup
					[
					color=cp3,
					style=dashed,
					]
					coordinates {
						(1,413.9320)
						(2,206.9660)
						(4,103.4830)
						(8,51.7415)
						(16,25.8708)
						(32,12.9354)
						(64,6.4677)
						(128,3.2338)
						(256,1.6169)
						(512,0.8085)
						(1024,0.4042)
					};
					\addlegendentry{PIPS-IPM++}
					\addlegendentry{Gurobi}
					\addlegendentry{PaPILO}
					\addlegendentry{Ideal linear}
				\end{loglogaxis}
			\end{tikzpicture}
		\end{adjustbox}
	\end{subfigure}%
	\begin{subfigure}{.5\textwidth}
		\caption{simple\_tudpaper1}
		\begin{adjustbox}{width=0.95\textwidth}
			\begin{tikzpicture}
				\begin{loglogaxis}
					[
					xlabel={\MPI processes/threads},
					ylabel={Presolve time [s]},
					log ticks with fixed point,
					scale only axis,
					xmin=1, xmax=1024,
					ymin=0.1, ymax=1024,
					legend pos=south west,
					grid style={line width=1pt, draw=black, opacity=0.1},
					xmajorgrids={true},
					ymajorgrids={true},
					log basis x=2
					]
					\addplot %beam-me-mext/simple_instances/SIMPLE_tudpaper1/allblocksPips
					[
					color=cp6,
					mark=square*,
					thick
					]
					coordinates {
						(1,383.390)
						(2,188.953)
						(4,103.066)
						(8,52.470)
						(16,27.311)
						(32,19.430)
						(64,12.823)
						(128,15.024)
						(256,14.568)
						(438,10.673)
					};
					\addplot %gurobi
					[
					color=cp5,
					mark=triangle*,
					thick
					]
					coordinates {
						(1,227.04)
						(2,227.55)
						(4,224.23)
						(8,225.78)
						(16,224.80)
						(32,226.66)
						(64,225.65)
						(128,224.84)
					};
					\addplot % papilo
					[
					color=cp2,
					mark=diamond*,
					thick
					]
					coordinates {
						(1,183.500)
						(2,176.889)
						(4,136.187)
						(8,137.913)
						(16,138.831)
						(32,110.705)
						(64,138.246)
						(128,117.608)
					};
					\addplot % ideal speedup
					[
					color=cp3,
					style=dashed,
					]
					coordinates {
						(1,383.3900)
						(2,191.6950)
						(4,95.8475)
						(8,47.9237)
						(16,23.9619)
						(32,11.9809)
						(64,5.9905)
						(128,2.9952)
						(256,1.4976)
						(438,0.8753)
					};
					\addlegendentry{PIPS-IPM++}
					\addlegendentry{Gurobi}
					\addlegendentry{PaPILO}
					\addlegendentry{Ideal linear}
				\end{loglogaxis}
			\end{tikzpicture}
		\end{adjustbox}
	\end{subfigure}
	\begin{subfigure}{.5\textwidth}
		\caption{elmod\_cwe15}
		\begin{adjustbox}{width=0.95\textwidth}
			\begin{tikzpicture}
				\begin{loglogaxis}
					[
					xlabel={\MPI processes/threads},
					ylabel={Presolve time [s]},
					log ticks with fixed point,
					scale only axis,
					xmin=1, xmax=1024,
					ymin=0.1, ymax=1024,
					legend pos=south west,
					grid style={line width=1pt, draw=black, opacity=0.1},
					xmajorgrids={true},
					ymajorgrids={true},
					log basis x=2
					]
					\addplot %beam-me-mext/TUD/Set4_CWE15/ELMOD_438_20_noVEnames
					[
					color=cp6,
					mark=square*,
					thick
					]
					coordinates {
						(1,428.812)
						(2,144.897)
						(4,112.297)
						(8,61.444)
						(16,32.100)
						(32,17.890)
						(64,9.896)
						(128,11.728)
						(256,5.925)
						(438,4.344)
					};
					\addplot %gurobi
					[
					color=cp5,
					mark=triangle*,
					thick
					]
					coordinates {
						(1,342.58)
						(2,341.36)
						(4,338.43)
						(8,340.98)
						(16,339.95)
						(32,341.90)
						(64,341.08)
						(128,341.47)
					};
					\addplot % papilo
					[
					color=cp2,
					mark=diamond*,
					thick
					]
					coordinates {
						(1,511.628)
						(2,477.853)
						(4,427.598)
						(8,411.380)
						(16,397.282)
						(32,403.112)
						(64,392.289)
						(128,418.107)
					};
					\addplot % ideal speedup
					[
					color=cp3,
					style=dashed,
					]
					coordinates {
						(1,428.8120)
						(2,214.4060)
						(4,107.2030)
						(8,53.6015)
						(16,26.8008)
						(32,13.4004)
						(64,6.7002)
						(128,3.3501)
						(256,1.6750)
						(438,0.9790)
					};
					\addlegendentry{PIPS-IPM++}
					\addlegendentry{Gurobi}
					\addlegendentry{PaPILO}
					\addlegendentry{Ideal linear}
				\end{loglogaxis}
			\end{tikzpicture}
		\end{adjustbox}
	\end{subfigure}%
	\begin{subfigure}{.5\textwidth}
		\caption{elmod\_set5}
		\begin{adjustbox}{width=0.95\textwidth}
			\begin{tikzpicture}
				\begin{loglogaxis}
					[
					xlabel={\MPI processes/threads},
					ylabel={Presolve time [s]},
					log ticks with fixed point,
					scale only axis,
					xmin=1, xmax=1024,
					ymin=0.1, ymax=1024,
					legend pos=south west,
					grid style={line width=1pt, draw=black, opacity=0.1},
					xmajorgrids={true},
					ymajorgrids={true},
					log basis x=2
					]
					\addplot %beam-me-mext/TUD/Set5/ELMOD_438_20_noVEnames
					[
					color=cp6,
					mark=square*,
					thick
					]
					coordinates {
						(1,1318.160)
						(2,380.847)
						(4,294.670)
						(8,157.793)
						(16,80.349)
						(32,45.860)
						(64,26.187)
						(128,30.723)
						(256,15.701)
						(438,7.681)
					};
					\addplot %gurobi
					[
					color=cp5,
					mark=triangle*,
					thick
					]
					coordinates {
						(1,889.48)
						(2,885.94)
						(4,884.23)
						(8,882.67)
						(16,884.54)
						(32,887.08)
						(64,886.04)
						(128,882.70)
					};
					\addplot % ideal speedup
					[
					color=cp3,
					style=dashed,
					]
					coordinates {
						(1,1318.1600)
						(2,659.0800)
						(4,329.5400)
						(8,164.7700)
						(16,82.3850)
						(32,41.1925)
						(64,20.5963)
						(128,10.2981)
						(256,5.1491)
						(438,3.0095)
					};
					\addlegendentry{PIPS-IPM++}
					\addlegendentry{Gurobi}
					\addlegendentry{Ideal linear}
				\end{loglogaxis}
			\end{tikzpicture}
		\end{adjustbox}
	\end{subfigure}
	\begin{subfigure}{.5\textwidth}
		\caption{yssp\_exp\_dcopf\_488r\_96b}
		\begin{adjustbox}{width=0.95\textwidth}
			\begin{tikzpicture}
				\begin{loglogaxis}
					[
					xlabel={\MPI processes/threads},
					ylabel={Presolve time [s]},
					log ticks with fixed point,
					scale only axis,
					xmin=1, xmax=1024,
					ymin=0.1, ymax=1024,
					legend pos=south west,
					grid style={line width=1pt, draw=black, opacity=0.1},
					xmajorgrids={true},
					ymajorgrids={true},
					log basis x=2
					]
					\addplot %beam-me-mext/yssp/YSSP_exp_DCopf_TaW_96b/blocks_YSSP_exp_DCopf_TaW_488r_1h_t_96b_novenames
					[
					color=cp6,
					mark=square*,
					thick
					]
					coordinates {
						(1,638.795)
						(2,300.482)
						(4,163.539)
						(8,84.492)
						(16,43.821)
						(32,24.886)
						(64,17.080)
						(96,16.331)
					};
					\addplot %gurobi
					[
					color=cp5,
					mark=triangle*,
					thick
					]
					coordinates {
						(1,323.52)
						(2,322.55)
						(4,321.16)
						(8,322.16)
						(16,322.63)
						(32,324.38)
						(64,321.84)
						(128,323.97)
					};
					\addplot % ideal speedup
					[
					color=cp3,
					style=dashed,
					]
					coordinates {
						(1,638.7950)
						(2,319.3975)
						(4,159.6987)
						(8,79.8494)
						(16,39.9247)
						(32,19.9623)
						(64,9.9812)
						(96,6.6541)
					};
					\addlegendentry{PIPS-IPM++}
					\addlegendentry{Gurobi}
					\addlegendentry{Ideal linear}
				\end{loglogaxis}
			\end{tikzpicture}
		\end{adjustbox}
	\end{subfigure}%
	\begin{subfigure}{.5\textwidth}
		\caption{yssp\_exp\_488r\_250b}
		\begin{adjustbox}{width=0.95\textwidth}
			\begin{tikzpicture}
				\begin{loglogaxis}
					[
					xlabel={\MPI processes/threads},
					ylabel={Presolve time [s]},
					log ticks with fixed point,
					scale only axis,
					xmin=1, xmax=1024,
					ymin=0.1, ymax=1024,
					legend pos=south west,
					grid style={line width=1pt, draw=black, opacity=0.1},
					xmajorgrids={true},
					ymajorgrids={true},
					log basis x=2
					]
					\addplot %beam-me-mext/yssp/exp_488r_t_250b/YSSP_exp_488r_t_250b_novenames
					[
					color=cp6,
					mark=square*,
					thick
					]
					coordinates {
						(1,497.449)
						(2,264.731)
						(4,140.021)
						(8,71.182)
						(16,38.982)
						(32,23.312)
						(64,14.003)
						(128,18.149)
						(250,10.031)
					};
					\addplot %gurobi
					[
					color=cp5,
					mark=triangle*,
					thick
					]
					coordinates {
						(1,301.66)
						(2,299.34)
						(4,298.21)
						(8,300.81)
						(16,299.07)
						(32,304.51)
						(64,300.29)
						(128,298.90)
					};
					\addplot % ideal speedup
					[
					color=cp3,
					style=dashed,
					]
					coordinates {
						(1,497.4490)
						(2,248.7245)
						(4,124.3623)
						(8,62.1811)
						(16,31.0906)
						(32,15.5453)
						(64,7.7726)
						(128,3.8863)
						(250,1.9898)
					};
					\addlegendentry{PIPS-IPM++}
					\addlegendentry{Gurobi}
					\addlegendentry{Ideal linear}
				\end{loglogaxis}
			\end{tikzpicture}
		\end{adjustbox}
	\end{subfigure}
	\begin{subfigure}{.5\textwidth}
		\caption{yssp\_exp\_488r\_1460b}
		\begin{adjustbox}{width=0.95\textwidth}
			\begin{tikzpicture}
				\begin{loglogaxis}
					[
					xlabel={\MPI processes/threads},
					ylabel={Presolve time [s]},
					log ticks with fixed point,
					scale only axis,
					xmin=1, xmax=2048,
					ymin=0.1, ymax=4096,
					legend pos=south west,
					grid style={line width=1pt, draw=black, opacity=0.1},
					xmajorgrids={true},
					ymajorgrids={true},
					log basis x=2
					]
					\addplot %beam-me-mext/yssp/exp_488r_t_1460b/YSSP_exp_488r_t_1460b_novenames
					[
					color=cp6,
					mark=square*,
					thick
					]
					coordinates {
						(1,2973.120)
						(2,1640.650)
						(4,424.516)
						(8,310.554)
						(16,160.544)
						(32,295.150)
						(64,439.239)
						(128,577.888)
						(256,595.402)
						(512,597.844)
						(1024,604.972)
						(1460,588.406)
					};
					\addplot %gurobi
					[
					color=cp5,
					mark=triangle*,
					thick
					]
					coordinates {
						(1,300.74)
						(2,298.80)
						(4,299.60)
						(8,300.26)
						(16,298.96)
						(32,300.10)
						(64,300.44)
						(128,299.47)
					};
					\addplot % ideal speedup
					[
					color=cp3,
					style=dashed,
					]
					coordinates {
						(1,2973.1200)
						(2,1486.5600)
						(4,743.2800)
						(8,371.6400)
						(16,185.8200)
						(32,92.9100)
						(64,46.4550)
						(128,23.2275)
						(256,11.6137)
						(512,5.8069)
						(1024,2.9034)
						(1460,2.0364)
					};
					\addlegendentry{PIPS-IPM++}
					\addlegendentry{Gurobi}
					\addlegendentry{Ideal linear}
				\end{loglogaxis}
			\end{tikzpicture}
		\end{adjustbox}
	\end{subfigure}%
	\begin{subfigure}{.5\textwidth}
		\caption{eustem}
		\begin{adjustbox}{width=0.95\textwidth}
			\begin{tikzpicture}
				\begin{loglogaxis}
					[
					xlabel={\MPI processes/threads},
					ylabel={Presolve time [s]},
					log ticks with fixed point,
					scale only axis,
					xmin=1, xmax=1024,
					ymin=0.1, ymax=4096,
					legend pos=south west,
					grid style={line width=1pt, draw=black, opacity=0.1},
					xmajorgrids={true},
					ymajorgrids={true},
					log basis x=2
					]
					\addplot %beam-me-mext/PSI/Set3/eustem_672_11_8jac
					[
					color=cp6,
					mark=square*,
					thick
					]
					coordinates {
						(1,344.352)
						(2,216.559)
						(4,76.218)
						(8,75.144)
						(16,29.724)
						(32,29.648)
						(64,15.206)
						(70,14.861)
					};
					\addplot %gurobi
					[
					color=cp5,
					mark=triangle*,
					thick
					]
					coordinates {
						(1,216.12)
						(2,217.06)
						(4,215.99)
						(8,216.56)
						(16,217.74)
						(32,215.91)
						(64,215.03)
						(128,217.54)
					};
					\addplot % papilo
					[
					color=cp2,
					mark=diamond*,
					thick
					]
					coordinates {
						(1,1627.532)
						(2,1596.101)
						(4,1501.230)
						(8,1495.820)
						(16,1462.238)
						(32,1480.127)
						(64,1471.694)
						(128,1469.383)
					};
					\addplot % ideal speedup
					[
					color=cp3,
					style=dashed,
					]
					coordinates {
						(1,344.3520)
						(2,172.1760)
						(4,86.0880)
						(8,43.0440)
						(16,21.5220)
						(32,10.7610)
						(64,5.3805)
						(70,4.9193)
					};
					\addlegendentry{PIPS-IPM++}
					\addlegendentry{Gurobi}
					\addlegendentry{PaPILO}
					\addlegendentry{Ideal linear}
				\end{loglogaxis}
			\end{tikzpicture}
		\end{adjustbox}
	\end{subfigure}
	\caption{Presolve scaling behavior of PIPS-IPM++, PaPILO, and Gurobi.} \label{fig:presolve_scaling}
\end{figure}
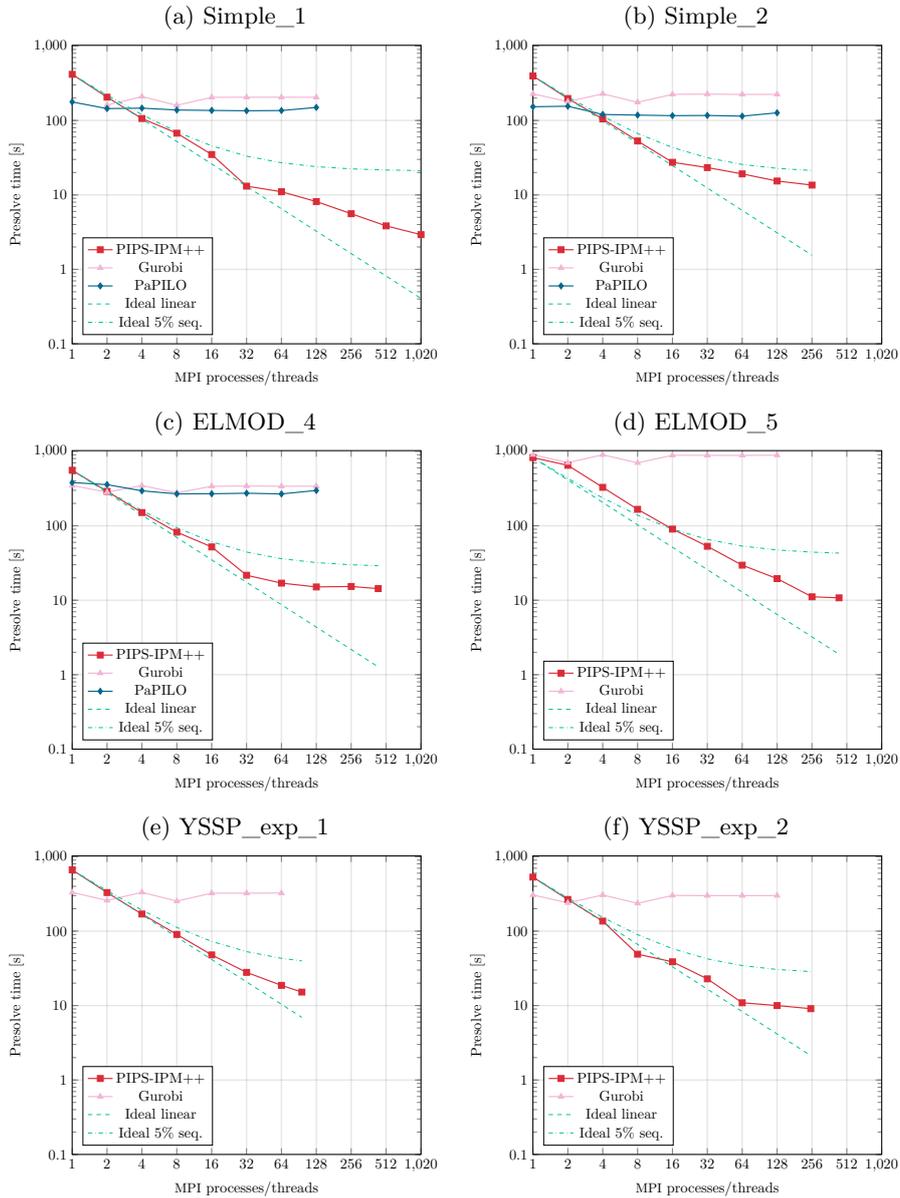
\fi

\subsection{Comparison against other presolvers}\label{sec:results_efficiency}

To evaluate the efficiency and efficacy of our implementation we compare the amount of problem reductions (removed nonzeros, constraints, and variables) and execution time of our presolve against the presolve of Gurobi and PaPILO. We ran PIPS-IPM++ using one node with 64 processes and with the maximum possible number of processes on multiple nodes, and PaPILO and Gurobi each with 64 threads on one node. For Gurobi and PaPILO, the sizes refer to the MPS file, which on 64 of the 82 instances has one additional row and column defining the objective function. Table~\ref{tab:compare_presolve_gur_pap_pips} shows the number of instances where presolve completes within the time limit, the presolve time, and the size of the reduced problem.

Independent runs of PIPS-IPM++ with the same number of processes produce identical presolved problems. Between 64 processes and the maximum number of processes, which distribute the blocks differently, the presolved problems differ on 8 of 81 instances, by at most 0.01\% of the nonzeros. Without enforcing the reduction order, two repeated runs gave identical presolved problems, while the presolve time decreased from 4.45\,s to 3.7\,s with 64 processes and from 4.06\,s to 3.55\,s and 4.02\,s with the maximum number of processes (sgm). This cost of determinism is concentrated on the five YSSP instances with 1460 blocks (219 to 494\,s vs.\ 37 to 87\,s with 64 processes). Without it, our presolve would be 7.3 instead of 6.1 times faster than Gurobi with 64 processes. On multiple nodes, the slowest of four runs took more than twice as long as the fastest on 14 of 71 instances, presumably due to network load and allocation topology. On a single node, times were stable. Since Gurobi and PaPILO are deterministic as well, we report the times of a single deterministic run.
\ifzibreport
\begin{table}[ht!]
\caption{Presolve results of PIPS-IPM++, Gurobi, and PaPILO. Time (sgm, shift 1\,s) and remaining problem size in percent of the original problem (sgm, shift 1\,\%) over the successfully presolved instances.}
\label{tab:compare_presolve_gur_pap_pips}
\adjustbox{max width=\textwidth}{%
\begin{tabular}{l|rrrrrr}
			\toprule
			\multicolumn{4}{c}{} & \multicolumn{3}{c}{reduced problem size} \\
			\cmidrule(lr){5-7} 
			Presolver  & \# success & \# failed & time (s) & nonzeros (\%) & constraints (\%) & variables (\%)\\
			\midrule
			PIPS-IPM++ (max proc.) & 82 & 0  &  4.15 & 68.43 & 63.07 & 77.47 \\
			PIPS-IPM++ (64 proc.)  & 82 & 0  &  4.39 & 68.43 & 63.07 & 77.47 \\
			Gurobi                 & 82 & 0  & 27.42 & 61.71 & 51.63 & 67.14 \\
			PaPILO                 & 70 & 12 & 82.54 & 66.71 & 61.90 & 72.60 \\
			\bottomrule
\end{tabular}
}
\end{table}
\else
\begin{table}[ht!]
\TABLE
{Presolve results of PIPS-IPM++, Gurobi, and PaPILO. Time (sgm, shift 1\,s) and remaining problem size in percent of the original problem (sgm, shift 1\,\%) over successfully presolved instances.\label{tab:compare_presolve_gur_pap_pips}}
{\begin{tabular}{l|rrrrrr}
\toprule
\multicolumn{4}{c}{} & \multicolumn{3}{c}{reduced problem size} \\
\cmidrule(lr){5-7} 
Presolver  & \# success & \# failed & time (s) & nonzeros (\%) & constraints (\%) & variables (\%)\\
\midrule
PIPS-IPM++ (max proc.) & 82 & 0  &  4.15 & 68.43 & 63.07 & 77.47 \\
PIPS-IPM++ (64 proc.)  & 82 & 0  &  4.39 & 68.43 & 63.07 & 77.47 \\
Gurobi     (64 proc.)  & 82 & 0  & 27.42 & 61.71 & 51.63 & 67.14 \\
PaPILO     (64 proc.)  & 70 & 12 & 82.54 & 66.71 & 61.90 & 72.60 \\
\bottomrule
\end{tabular}}
{}
\end{table}
\fi
While PaPILO fails on 12 instances, nine of which it cannot presolve within 1.5 hours, PIPS-IPM++ and Gurobi presolve all instances successfully. With 64 processes on one node, that is, on the same hardware as Gurobi and PaPILO, the presolve of PIPS-IPM++ is about a factor of six faster than Gurobi and a factor of 19 faster than PaPILO (sgm). On the 69 instances presolved by all three, the factors are 5.9 and 25. On the 18 instances with at least 100 million nonzeros, our presolve is 13 times faster than Gurobi with 64 processes and 21 times faster with the maximum number of processes (sgm, on the 17 instances presolved by both), while PaPILO only presolves 9 of these instances. Over the whole testset, the maximum number of processes reduces the sgm presolve time only slightly, from 4.39\,s to 4.15\,s: on 20 instances it is at most 64 anyway, and on most others presolve with 64 processes already takes only a few seconds, so that communication dominates with more processes. MPI communication takes about a fifth of the presolve time with up to 64 processes and about three quarters with more than 256 (\ifzibreport Appendix~\ref{appendix:mpi}\else online supplement, Appendix~D\fi).

Whether multiple compute nodes pay off depends mainly on the linking part of the problem rather than on its size or number of blocks. Over the 61 instances with more than 64 blocks, the speed-up from 64 processes to the maximum number of processes correlates only weakly with the number of nonzeros (Spearman correlation 0.43) and of blocks (0.24), but negatively with the share of linking constraints that couple more than two consecutive blocks and linking variables that appear in local constraints ($-0.65$). On the six instances of Section~\ref{sec:results_scaling} with small linking parts, the maximum number of processes is up to 3.4 times faster than 64 processes. Instances with comparatively large linking parts, such as the spde and remix\_nagsys\_de models, eustem, and the YSSP models with 1460 blocks, become slower: the global part of the problem is processed redundantly on every process, and the block contributions to the linking part (e.g., activities and bounds) are combined by communication whose cost grows with the number of processes.

Looking at the amount of nonzeros left in the problem after presolve, which is the most important metric for the performance of Simplex, IPM, and PDLP algorithms, PIPS-IPM++ reduces the problems to 68.4\% of their original nonzeros, PaPILO to 66.7\%, and Gurobi to 61.7\%. On the 69 instances presolved by all three (PaPILO fails on large YSSP instances on which PIPS-IPM++ removes many nonzeros), 71.4\%, 66.7\%, and 65.5\% remain, respectively. Gurobi also removes 11.4 percentage points more constraints and 10.3 more variables than PIPS-IPM++. This difference is not simply explained by the requirement to preserve the arrowhead structure, since the reductions of Gurobi and PaPILO largely preserve it as well (Section~\ref{sec:results_external}). Rather, both implement reductions that PIPS-IPM++ does not, such as the substitution of free variables in constraints with more than two entries, singleton stuffing, sparsification, dual inference, and dominated columns. As Section~\ref{sec:results_external} shows, these reductions do not make the IPM noticeably faster on the given instances.

Overall, we see a strong runtime advantage of PIPS-IPM++ over both the academic and the commercial alternatives: already on a single compute node, our presolve is about six times faster than Gurobi, and on large instances with small linking parts, multiple compute nodes increase this advantage to factors of 20 to 115. Gurobi's additional 6.7 percentage points of removed nonzeros thus cost about six times the presolve time on a single node, and they do not make the IPM noticeably faster (Section~\ref{sec:results_external}).

\subsection{Presolve impact in PIPS-IPM++}\label{sec:results_impact_ipm}

To investigate the impact of presolving on solving the models, we ran PIPS-IPM++ on all 82 instances of our testset with the maximum number of processes and a time limit of one hour, with and without presolve. Table~\ref{tab:pips_w_wo} shows the number of solved instances and, for the 35 instances solved in both runs, the time, the number of IPM iterations, and the time per iteration. The last column gives the sgm time over all instances, where each unsolved instance (time or iteration limit, wrong status, or crash) contributes 3600 seconds.
\ifzibreport
\begin{table}[ht!]
\begin{center}
\caption{PIPS-IPM++ performance with and without presolve. Times and iterations as sgm (shift 1\,s / 1); time per iteration: geometric mean; unsolved instances contribute 3600\,s to the last column.}
\label{tab:pips_w_wo}
\adjustbox{max width=\textwidth}{%
\begin{tabular}{l|rrrrr}
\toprule
\multicolumn{2}{c}{} & \multicolumn{3}{c}{solved by both (35)} & all (82) \\
\cmidrule(lr){3-5} \cmidrule(lr){6-6}
 & \# solved & time (s) & iterations & time/iteration (s) & time (s) \\
\midrule
Presolve     & 48 & 157.7 & 88.4 & 1.95 &  640.5 \\
No presolve  & 37 & 215.6 & 87.2 & 2.45 & 1068.6 \\
\bottomrule
\end{tabular}
}
\end{center}
\end{table}
\else
\begin{table}[ht!]
\TABLE
{PIPS-IPM++ performance with and without presolve. Times and iterations: sgm (shift 1\,s / 1); time per iteration: geometric mean; unsolved instances contribute 3600\,s to the last column.\label{tab:pips_w_wo}}
{\begin{tabular}{l|rrrrr}
\toprule
\multicolumn{2}{c}{} & \multicolumn{3}{c}{solved by both (35)} & all (82) \\
\cmidrule(lr){3-5} \cmidrule(lr){6-6}
 & \# solved & time (s) & iterations & time/iteration (s) & time (s) \\
\midrule
Presolve     & 48 & 157.7 & 88.4 & 1.95 &  640.5 \\
No presolve  & 37 & 215.6 & 87.2 & 2.45 & 1068.6 \\
\bottomrule
\end{tabular}}
{}
\end{table}
\fi
The impact of presolving for PIPS-IPM++ is significant. With presolve, PIPS-IPM++ solves 48 instances, without presolve 37. On the instances solved in both runs, presolve hardly changes the number of iterations, but each iteration is about 20\% faster on the presolved problems, and the time decreases by 27\%. Of the 13 instances solved only with presolve, the run without presolve reports infeasibility on three (wrongly, as they are solved with presolve), hits the iteration limit on three, the time limit on four, and crashes on three. Presolve thus also makes the IPM more robust. Two instances are solved only without presolve. With presolve, one hits the time limit and one the iteration limit. Over all instances, the sgm time decreases by 40\%. This highlights the importance of presolve within PIPS-IPM++. We also believe that other AHLP exploiting methods can benefit similarly from our contribution.

\subsection{Comparison with external presolve}\label{sec:results_external}

An alternative to our own presolve is to use a general-purpose presolver such as Gurobi's or PaPILO's and then pass the presolved model to PIPS-IPM++. Doing so, however, the modeler's annotation is lost, since neither Gurobi nor PaPILO provide a mapping from the presolved model to the original one. We evaluated two ways of annotating the presolved models of Section~\ref{sec:results_efficiency}. First, we annotated them automatically with pipstools (\url{https://gitlab.com/pips-ipmpp/pipstools}, commit \texttt{b381aae6}), an open-source tool that partitions the constraint matrix into the same number of blocks as the modeler's annotation using the hypergraph partitioner KaHyPar~\cite{Schlag2020_KaHyPar}. This is the option currently available to users of PIPS-IPM++ who cannot supply a manual annotation. We also re-annotated the original models using pipstools. Second, we transferred the modeler's annotation to the presolved models: each variable is assigned to the block of its original counterpart, or becomes a linking variable if there is none, and each constraint is positioned according to the blocks of its variables. We reconstructed the correspondence from the variable names, which both presolvers preserve, and verified it by the nonzero counts of all variables. This inherited annotation is fragile: it relies on the presolver keeping the variable names, it cannot assign variables introduced by the presolver, and it requires the annotation of the complete model in a single file. We could use it for 62 instances, all presolved by Gurobi, of which PaPILO presolved 56 within its time limit. For the other 20 instances, the annotation exists only as split block files without names. We solved all models with PIPS-IPM++, with and without our presolve, as in Section~\ref{sec:results_impact_ipm}.

The performance of PIPS-IPM++ depends critically on the share of linking constraints that couple only two consecutive blocks (2-link constraints), which PIPS-IPM++ treats in a block-sprase fashion, while all other linking constraints enter the Schur complement in as dense sub-blocks~\cite{kempke2024massivelyparallelinteriorpointmethodarrowhead}. Table~\ref{tab:external} shows, for each combination of model and annotation, the number of solved instances together with the structure of the annotation. The reference columns compare each row directly with our presolve.
\ifzibreport
\begin{table}[ht!]
\centering
\caption{PIPS-IPM++ on externally presolved models. Solved instances without and with our presolve, compared with the modeler's annotation of the original models on the same instances (reference). Linking constraints and 2-link share of the annotation as read by PIPS-IPM++ (medians).}
\label{tab:external}
\adjustbox{max width=\textwidth}{%
\begin{tabular}{ll|r|rr|rr|rr}%r}
\toprule
 & & & \multicolumn{2}{c|}{\# solved} & \multicolumn{2}{c|}{reference} & linking & 2-link \\% & no 2-link \\
annotation & presolver & \# inst. & w/o pres. & w/ pres. & w/o pres. & w/ pres. & constraints & share \\%& treatment \\
\midrule
modeler   & --               & 82 & 37 & 48 & 37 & 48 &  17{,}694 & 0.99 \\% & 12 \\
automatic & --               & 81 & 13 & 21 & 37 & 48 &  95{,}505 & 0.48 \\%& 55 \\
automatic & Gurobi           & 77 & 13 & 16 & 34 & 45 & 102{,}156 & 0.46 \\%& 54 \\
automatic & PaPILO           & 69 & 13 & 14 & 35 & 46 & 109{,}426 & 0.46 \\%& 50 \\
inherited & Gurobi           & 62 & 40 & 43 & 32 & 39 &  21{,}351 & 1.00 \\%&  8 \\
inherited & PaPILO           & 56 & 40 & 37 & 32 & 39 &  10{,}251 & 1.00 \\%& 12 \\
\bottomrule
\end{tabular}
}
\end{table}
\else
\begin{table}[ht!]
\TABLE
{PIPS-IPM++ on externally presolved models. Solved instances without and with our presolve, compared with the modeler's annotation of the original models on the same instances (reference). Linking constraints and 2-link share of the annotation as read by PIPS-IPM++ (medians). \label{tab:external}}
{\begin{tabular}{ll|r|rr|rr|rr}%r}
\toprule
 & & & \multicolumn{2}{c|}{\# solved} & \multicolumn{2}{c|}{reference} & linking & 2-link \\%& no 2-link \\
annotation & presolver & \# inst. & w/o pres. & w/ pres. & w/o pres. & w/ pres. & constraints & share \\% & treatment \\
\midrule
modeler   & --               & 82 & 37 & 48 & 37 & 48 &  17{,}694 & 0.99 \\%& 12 \\
automatic & --               & 81 & 13 & 21 & 37 & 48 &  95{,}505 & 0.48 \\%& 55 \\
automatic & Gurobi           & 77 & 13 & 16 & 34 & 45 & 102{,}156 & 0.46 \\%& 54 \\
automatic & PaPILO           & 69 & 13 & 14 & 35 & 46 & 109{,}426 & 0.46 \\%& 50 \\
inherited & Gurobi           & 62 & 40 & 43 & 32 & 39 &  21{,}351 & 1.00 \\%&  8 \\
inherited & PaPILO           & 56 & 40 & 37 & 32 & 39 &  10{,}251 & 1.00 \\%& 12 \\
\bottomrule
\end{tabular}}
{}
\end{table}
\fi
With the automatic annotation, PIPS-IPM++ solves at most 21 instances, compared to 48 instances with the modeler's annotation and our presolve. The automatic annotation detects about six times as many linking constraints, only about half of them 2-link constraints (99\% for the modeler's annotation), and most of the failed runs end with an out-of-memory error. Since the automatic annotation of the original models yields the same structure, the loss is caused by the re-annotation that an external presolve requires and that our presolve avoids.

With the inherited annotation, the externally presolved models mostly keep the arrowhead structure (Table~\ref{tab:external}). This is a best case that is not available to users and that we could construct for only 62 of the 82 instances. Even so, the external presolves solve about as many instances as our presolve on the same instances, with similar times of PIPS-IPM++ (Table~\ref{tab:inherited}), and including the time of the external presolve, our presolve is the fastest (124\,s vs.\ 132\,s for Gurobi and 141\,s vs.\ 219\,s for PaPILO). The additional reductions done by Gurobi do speed-up the IPM (124\,s down to 118\,s), but are absorbed by longer presolve times and in addition, can damage the structure: on two ELMOD instances, Gurobi creates about 680{,}000 variables that belong to no block and become linking variables. Without our presolve, PIPS-IPM++ runs out of memory on these instances, with it, PIPS-IPM++ takes about 3.7 times as long as with our presolve alone (3064\,s and 3259\,s vs.\ 819\,s and 894\,s).
\ifzibreport
\begin{table}[ht!]
\centering
\caption{PIPS-IPM++ with the inherited annotation: solved instances and time (sgm) on the instances solved by all four presolve combinations, without and with the time of the external presolve.}
\label{tab:inherited}
\begin{tabular}{l|rrr|rrr}
\toprule
 & \multicolumn{3}{c|}{Gurobi (62 inst., 28 common)} & \multicolumn{3}{c}{PaPILO (56 inst., 30 common)} \\
presolve & \# solved & PIPS-IPM++ & total & \# solved & PIPS-IPM++ & total \\
\midrule
none            & 32 & 137.9 & 137.9 & 32 & 166.4 & 166.4 \\
ours            & 39 & 124.2 & 124.2 & 39 & 141.4 & 141.4 \\
external        & 40 & 117.7 & 132.1 & 40 & 142.1 & 219.1 \\
external + ours & 43 & 117.9 & 131.9 & 37 & 149.7 & 222.2 \\
\bottomrule
\end{tabular}
\end{table}
\else
\begin{table}[ht!]
\TABLE
{PIPS-IPM++ with the inherited annotation: solved instances and time (sgm) on the instances solved by all four presolve combinations, without and with the time of the external presolve. \label{tab:inherited}}
{\begin{tabular}{l|rrr|rrr}
\toprule
 & \multicolumn{3}{c|}{Gurobi (62 inst., 28 common)} & \multicolumn{3}{c}{PaPILO (56 inst., 30 common)} \\
presolve & \# solved & PIPS-IPM++ & total & \# solved & PIPS-IPM++ & total \\
\midrule
none            & 32 & 137.9 & 137.9 & 32 & 166.4 & 166.4 \\
ours            & 39 & 124.2 & 124.2 & 39 & 141.4 & 141.4 \\
external        & 40 & 117.7 & 132.1 & 40 & 142.1 & 219.1 \\
external + ours & 43 & 117.9 & 131.9 & 37 & 149.7 & 222.2 \\
\bottomrule
\end{tabular}}
{}
\end{table}
\fi
In our experiments, external presolve's cost grows with the size of the model. On the 18 instances with at least 100 million nonzeros, Gurobi needs 94\,s (in sgm, at most 308\,s) and PaPILO 504\,s only to read the MPS file, while PIPS-IPM++ reads the model from its block files in 4.6\,s with 64 processes and in 2.0\,s with the maximum number of processes. Gurobi's presolve then takes 247\,s on the 17 instances it presolves, compared to 19.2\,s and 11.8\,s for our presolve, and PaPILO presolves only 9 of the 18 instances. On the ELMOD instances, Gurobi's presolve takes ~40\% of the total solution time of PIPS-IPM++ with our presolve.

Most importantly, a sequential presolve breaks the distributed pipeline of model generation, reading, and solution which PIPS-IPM++ was designed for: the complete model has to be generated all at once and written to a single file, and the presolved model has to be re-annotated and split into blocks again. This can be avoided by using the distributed model generation available with PIPS-IPM++ GAMS interface~\cite{Scholz2020_BestPracticeGuide}. A general-purpose presolve that carries the annotation through its reductions would avoid the re-annotation, but not these limitations, and still nothing prevents it from introducing variables that belong to no block, as on the ELMOD instances.

\section{Conclusion}\label{sec:conclusion}

We presented and evaluated our distributed parallel, structure preserving presolve implemented in the parallel solver PIPS-IPM++. We demonstrated its scalability and general effectiveness on a large set of AHLPs showing that already with 64 processes on a single compute node our presolve outperforms Gurobi by a factor of six and PaPILO by a factor of 19, while removing two to five percentage points fewer nonzeros than PaPILO. On large instances with small linking parts, deploying our presolve on multiple compute nodes makes it 20 to 115 times faster than Gurobi's presolve. Presolve also makes PIPS-IPM++ more robust: it solves 48 instead of 37 instances of our testset. Compared with an external presolve by Gurobi or PaPILO, our presolve does not require recovering the block structure, which an automatic annotation cannot do well. Even with the modeler's annotation carried over to the externally presolved models, PIPS-IPM++ is not faster overall, the external presolve is slow on large models, and it breaks the distributed pipeline of model generation, reading, and solution, which our presolve keeps intact. We demonstrated that structure specific implementations of algorithms can help speed-up the solution process.

For our presolve implementation we mainly see two ways forward. First, make the implementation easily accessible to other structure exploiting methods that work on block structured LPs. Second, improve presolve further by adding more presolvers and tuning. Our experiments comparing against Gurobi indicate a gap in algorithmic presolve that we believe is not simply explained by block-structure preservation restrictions.

%\THEEndNotes
%\begingroup \parindent 0pt \parskip 0.0ex \def\enotesize{\normalsize} \theendnotes \endgroup

% Appendix here
% Options are (1) APPENDIX (with or without general title) or
%             (2) APPENDICES (if it has more than one unrelated sections)
% Outcomment the appropriate case if necessary
%

\ifzibreport
\appendix
\else

\fi
%

% Acknowledgments here
\ACKNOWLEDGMENT{The work for this article has been conducted in the Research Campus MODAL funded by the German Federal Ministry of Research, Technology and Space (BMFTR) (fund numbers 05M14ZAM, 05M20ZBM, 05M2025). Further, the described research activities were funded by the German Federal Ministry for Economic Affairs and Energy (BMWi) within the project PEREGRINE (grant number 
03EI1082B). The authors gratefully acknowledge the computational and data resources provided by the Leibniz Supercomputing Centre (www.lrz.de). The authors used the generative AI tool Claude Code (Anthropic Opus 5.5) to write evaluation, job, and packaging scripts and for language editing. The authors checked all scripts, results, and text and take full responsibility for the content.}

% References here (outcomment the appropriate case)

% CASE 1: BiBTeX used to constantly update the references
%   (while the paper is being written).
\bibliographystyle{IJOC_style/informs2014.bst} % outcomment this and next line in Case 1
\bibliography{references.bib} % if more than one, comma separated

% CASE 2: BiBTeX used to generate mypaper.bbl (to be further fine tuned)
%\input{mypaper.bbl} % outcomment this line in Case 2

%\bibliographystyle{nonumber}

%%%%%%%%%%%%%%%%%
\end{document}